\documentclass{article}
\usepackage{amsmath,amssymb}
\newtheorem{df}{Definition}[section]
\newtheorem{thm}[df]{Theorem}
\newtheorem{lem}[df]{Lemma}
\newtheorem{prop}[df]{Proposition}
\newtheorem{cor}[df]{Corollary}
\newtheorem{exa}[df]{Example}
\newtheorem{con}[df]{Conjecture}
\title{Directional wavelets on $n$--dimensional spheres}
\author{Ilona Iglewska-Nowak\footnote{West Pomeranian University of Technology, School of Mathematics, al. Piast\'ow 17, 70--310 Szczecin, Poland}}

\begin{document}

\maketitle

\bibliographystyle{amsplain}

\begin{abstract}
Directional Poisson wavelets, being directional derivatives of Poisson kernel, are introduced on $n$--dimensional spheres. It is shown that, slightly modified and together with another wavelet family, they are an admissible wavelet pair according to the definition derived from the theory of approximate identities. We investigate some of the properties of directional Poisson wavelets, such as recursive formulae for their Fourier coefficients or explicit representations as functions of spherical variables (for some of the wavelets). We derive also an explicit formula for their Euclidean limits.
\end{abstract}

\begin{bfseries}Keywords:\end{bfseries}spherical directional wavelets, Poisson multipole wavelets, $n$-spheres \\
\begin{bfseries}AMS Classification:\end{bfseries} 42C40, 42B20

\section{Introduction}

In the present paper we continue our investigation of wavelets over $n$--dimensional spheres. The definition we base on comes from~\cite{EBCK09} and~\cite{IIN14CWT} and it is a generalization of definitions derived from the theory of approximate identities and singular integrals for the two--dimensional case \cite{FW,FW-C,FGS-book} as well as for the three-- and $n$--dimensional cases \cite{sB09,BE10WS3,BE10KBW,sE11} .

Inspired by~\cite{HH09} we present nonzonal wavelets being directional derivatives of Poisson kernel. Nonzonal spherical wavelets derived from approximate identities were introduced 2009~\cite{EBCK09} and to our knowledge the present research is the first attempt to define a concrete wavelet family satisfying the conditions of the definition given in that paper. The motivation for the choice of Poisson wavelets are their excellent properties in the zonal case (zonal Poisson wavelets are derivatives of Poisson kernel along the axis through the origin of the sphere and the origin of the Poisson kernel \cite{HI07}): they possess explicit representations and are well--localized \cite{IIN14PW}, they also have discrete frames \cite{IH10,IIN14WF} and therefore, are well-suited for computations \cite{HCM03,CPMDHJ05}. Also directional wavelets reveal some of that properties, as it is shown in the present paper: the kernel of the wavelet transform is a linear combination of Mexican needlets~\cite{aM10} (and the Gauss--Weierstrass wavelet), thus, it is very well localized; it is straightforward to obtain an explicit representation of the first directional wavelet, cf. Example~\ref{exa:g_rho_1}, a computation for the second--order wavelet is contained in the Appendix. In a similar way explicit representations for higher--order wavelets can be obtained. An existence of discrete frames needs to be investigated, but the author supposes it is provided due to the good localization of the kernel of the wavelet transform (similarly in the case of zonal Poisson wavelets, the existence proof for discrete frames is based on the kernel localization~\cite{IIN14WF}). Further, the Euclidean limits of directional Poisson wavelets exists and is given by a simple formula, cf. Section~\ref{sec:Euclidean_limit}.

The paper is organized as follows. Section~\ref{sec:sphere} contains basic information about analysis of functions on spheres. We introduce directional Poisson wavelets as directional derivatives of Poisson kernel in Section~\ref{sec:definition} and compute recursive formulae for their series representation in Section~\ref{sec:derivatives_of_rotations}. It is shown in Section~\ref{sec:wavelets} that certain linear combinations of that Poisson wavelets satisfy conditions of the definition of wavelets derived from an approximate identity. Finally, Euclidean limit of directional Poisson wavelets are computed in Section~\ref{sec:Euclidean_limit}. It is shown in Appendix how to obtain an explicit representation of a directional Poisson wavelet on the example of the second directional derivative of Poisson kernel.

\section{Preliminaries}\label{sec:sphere}

\subsection{Functions on the sphere}

By $\mathcal{S}^n$ we denote the $n$--dimensional unit sphere in $n+1$--dimensional Euclidean space~$\mathbb{R}^{n+1}$ with the rotation--invariant measure~$d\sigma$ normalized such that
$$
\Sigma_n=\int_{\mathcal{S}^n}d\sigma=\frac{2\pi^{(n+1)/2}}{\Gamma\bigl((n+1)/2\bigr)}.
$$
The surface element $d\sigma$ is explicitly given by
$$
d\sigma=\sin^{n-1}\theta_1\,\sin^{n-2}\theta_2\dots\sin\theta_{n-1}d\theta_1\,d\theta_2\dots d\theta_{n-1}d\varphi,
$$
where $(\theta_1,\theta_2,\dots,\theta_{n-1},\varphi)\in[0,\pi]^{n-1}\times[0,2\pi)$ are spherical coordinates satisfying
\begin{equation}\label{eq:spherical_coordinates}\begin{split}
x_1&=\cos\theta_1,\\
x_2&=\sin\theta_1\cos\theta_2,\\
x_3&=\sin\theta_1\sin\theta_2\cos\theta_3,\\
&\dots\\
x_{n-1}&=\sin\theta_1\sin\theta_2\dots\sin\theta_{n-2}\cos\theta_{n-1},\\
x_n&=\sin\theta_1\sin\theta_2\dots\sin\theta_{n-2}\sin\theta_{n-1}\cos\varphi,\\
x_{n+1}&=\sin\theta_1\sin\theta_2\dots\sin\theta_{n-2}\sin\theta_{n-1}\sin\varphi.
\end{split}\end{equation}
$\left<x,y\right>$ or $x\cdot y$ stand for the scalar product of vectors with origin in~$O$ and endpoint on the sphere. As long as it does not lead to misunderstandings, we identify these vectors with points on the sphere.

The $\mathcal{L}^p(\mathcal{S}^n)$--norm of a function is given by
$$
\|f\|_{\mathcal{L}^p(\mathcal{S}^n)}=\left[\frac{1}{\Sigma_n}\int_{\mathcal S^n}|f(x)|^p\,d\sigma(x)\right]^{1/p}.
$$
The scalar product of $f,g\in\mathcal L^2(\mathcal S^n)$ is defined by
$$
\left<f,g\right>_{\mathcal L^2(\mathcal S^n)}=\frac{1}{\Sigma_n}\int_{\mathcal S^n}\overline{f(x)}\,g(x)\,d\sigma(x),
$$
such that $\|f\|_2^2=\left<f,f\right>$.

A function is called zonal if its value depends only on $\theta=\theta_1=\left<\widehat e,x\right>$, where~$\widehat e$ is the north-pole of the sphere
$$
\widehat e=(1,0,0,\dots,0).
$$
It is invariant with respect to the rotation about the axis through~$O$ and~$\widehat e$. The subspace of $p$--integrable zonal functions is isomorphic to and will be identified with the space~$\mathcal L_\lambda^p$ of functions over the interval $[-1,1]$ having finite norm
$$
\|f\|_{\mathcal{L}_\lambda^p([-1,1])}=\left[\frac{\Sigma_{2\lambda}}{\Sigma_{2\lambda+1}}\int_{-1}^1|f(t)|^p\left(1-t^2\right)^{\lambda-1/2}dt\right]^{1/p}.
$$
The norms satisfy
$$
\|f\|_{\mathcal{L}^p(\mathcal{S}^n)}=\|f\|_{\mathcal{L}_\lambda^p([-1,1])},
$$
where~$\lambda$ and~$n$ are related by
$$
\lambda=\frac{n-1}{2}.
$$
We identify zonal functions with functions over the interval $[-1,1]$, i.e., whenever it does not lead to misunderstandings, we write
$$
f(x)=f(\cos\theta_1).
$$

Gegenbauer polynomials $C_l^\lambda$ of order~$\lambda\in\mathbb R$ and degree $l\in\mathbb{N}_0$, are defined in terms of their generating function
\begin{equation}\label{eq:generating_fct}
\sum_{l=0}^\infty C_l^\lambda(t)\,r^l=\frac{1}{(1-2tr+r^2)^\lambda},\qquad t\in[-1,1],
\end{equation}
they also satisfy the relation
$$
\sum_{l=0}^\infty\frac{\lambda+l}{\lambda}\,C_l^\lambda(t)\,r^l=\frac{1-r^2}{(1-2tr+r^2)^{\lambda+1}},\qquad t\in[-1,1].
$$
Explicitly we have
\begin{equation}\label{eq:C_jako_szereg}
C_l^\lambda(t)=\sum_{j=0}^{[l/2]}(-1)^j\frac{\Gamma(\lambda+l-j)}{j!\,(l-2j)!\,\Gamma(\lambda)}\,(2t)^{l-2j},
\end{equation}
cf. \cite[Sec. IX.3.1, formula (3)]{Vilenkin}.
A set of Gegenbauer polynomials $\bigl\{C_l^\lambda\bigr\}_{l\in\mathbb N_0}$ builds a complete orthogonal system on $[-1,1]$ with weight $(1-t^2)^{\lambda-1/2}$. Consequently, it is an orthogonal basis for zonal functions on the $(2\lambda+1)$--dimensional sphere. The following relations are valid for Gegenbauer polynomials:
\begin{align}
&(l+1)\,C_{l+1}^\lambda(t)=2(\lambda+l)\,t\,C_l^\lambda(t)-(2\lambda+l-1)\,C_{l-1}^\lambda(t)\label{eq:recGg1}\\
&(l+1)\,C_{l+1}^{\lambda-1}(t)=2\,(\lambda-1)\left[t\,C_l^\lambda(t)-C_{l-1}^\lambda(t)\right]\label{eq:recGg2}\\
&l\,C_l^\lambda(t)=(2\lambda+l-1)\,t\,C_{l-1}^\lambda(t)-2\lambda\,(1-t^2)\,C_{l-2}^{\lambda+1}(t)\label{eq:recGg3}\\
&\frac{d}{dt}\,C_l^\lambda(t)=2\lambda\,C_{l-1}^{\lambda+1}(t)\label{eq:derGg}
\end{align}
cf. \cite{GR}, formulae~8.933.1--2,4 and~8.939.1.

Let $Q_l$ denote a polynomial on~$\mathbb{R}^{n+1}$ homogeneous of degree~$l$, i.e., such that $Q_l(az)=a^lQ_l(z)$ for all $a\in\mathbb R$ and $z\in\mathbb R^{n+1}$, and harmonic in~$\mathbb{R}^{n+1}$, i.e., satisfying $\nabla^2Q_l(z)=0$, then $Y_l(x)=Q_l(x)$, $x\in\mathcal S^n$, is called a hyperspherical harmonic of degree~$l$. The set of hyperspherical harmonics of degree~$l$ restricted to~$\mathcal S^n$ is denoted by $\mathcal H_l(\mathcal S^n)$. Hyperspherical harmonics of distinct degrees are orthogonal to each other. The number of linearly independent hyperspherical harmonics of degree~$l$ is equal to
$$
N=N(n,l)=\frac{(n+2l-1)(n+l-2)!}{(n-1)!\,l!}.
$$
Addition theorem states that
$$
C_l^\lambda(x\cdot y)=\frac{\lambda}{\lambda+l}\,\sum_{\kappa=1}^N\overline{Y_l^\kappa(x)}\,Y_l^\kappa(y)
$$
for any orthonormal set $\{Y_l^\kappa\}_{\kappa=1,2,\dots,N(n,l)}$ of hyperspherical harmonics of degree~$l$ on~$\mathcal S^n$, cf. discussion in~\cite{IIN14CWT}.

In this paper, we will be working with the orthonormal basis for~$\mathcal L^2(\mathcal S^n)=\overline{\bigoplus_{l=0}^\infty\mathcal H_l}$, consisting of hyperspherical harmonics given by
\begin{equation}\label{eq:spherical_harmonics}
Y_l^k(x)=A_l^k\prod_{\tau=1}^{n-1}C_{k_{\tau-1}-k_\tau}^{\frac{n-\tau}{2}+k_\tau}(\cos\theta_\tau)\sin^{k_\tau}\!\theta_\tau\cdot e^{\pm ik_{n-1}\varphi}
\end{equation}
with $l=k_0\geq k_1\geq\dots\geq k_{n-1}\geq0$, $k$ being a sequence $(k_1,\dots,\pm k_{n-1})$ of integer numbers, and normalization constants
$$
A_l^k=\left(\frac{1}{\Gamma\left(\frac{n+1}{2}\right)}\prod_{\tau=1}^{n-1}\frac{2^{n-\tau+2k_\tau-2}\,(k_{\tau-1}-k_\tau)!\,(n-\tau+2k_{\tau-1})\,\Gamma^2(\frac{n-\tau}{2}+k_\tau)}{\sqrt\pi\,(n-\tau+k_{\tau-1}+k_\tau-1)!}\right)^{1/2},
$$
compare~\cite[Sec. IX.3.6, formulae (4) and (5)]{Vilenkin}. The set of non-increasing sequences~$k$ in $\mathbb N_0^{n-1}\times\mathbb Z$ with elements bounded by~$l$ will be denoted by $\mathcal M_{n-1}(l)$. For $k=(k_1,0,\dots,0)$ we have by doubling formula for gamma function \cite[8.335.1]{GR}
\begin{equation}\label{eq:Alk1}
A_l^k=\sqrt{\frac{2^{2n+2k_1-6}\,(l-k_1)!\,k_1!\,(n+2l-1)\,(n+2k_1-2)\,\Gamma^2(\frac{n-1}{2}+k_1)\,\Gamma^2(\frac{n-2}{2})}
   {(n-1)\,\pi\,(n+l+k_1-2)!\,(n+k_1-3)!}}.
\end{equation}
In particular, for $k=(0,0,\dots,0)$,
\begin{equation}\label{eq:Al0}
A_l^0=\sqrt{\frac{(n-2)!\,l!\,(n+2l-1)}{(n+l-2)!\,(n-1)}},
\end{equation}
cf. \cite[Sec. IX.3.6, formula (6)]{Vilenkin}. In the two--dimensional case, formula~\eqref{eq:Alk1} simplifies to
\begin{equation}\label{eq:Alk_n2}
A_l^k=2^k\Gamma\left(k+\frac{1}{2}\right)\sqrt{\frac{(2l+1)(l-k)!}{\pi\,(l+k)!}}.
\end{equation}

Further, Funk--Hecke theorem states that for $f\in\mathcal L_\lambda^1\bigl([-1,1])$ and $Y_l\in\mathcal H_l(\mathcal S^n)$, $l\in\mathbb{N}_0$,
\begin{equation}\label{eq:FH}\begin{split}
\int_{\mathcal S^n}&Y_l(y)\,f(x\cdot y)\,d\sigma(y)\\
&=Y_l(x)\cdot\frac{(4\pi)^\lambda\Gamma(l+1)\Gamma(\lambda)}{\Gamma(2\lambda+l)}\int_{-1}^1 f(t)\,C_l^\lambda(t)\left(1-t^2\right)^{\lambda-1/2}dt.
\end{split}\end{equation}

Every $\mathcal{L}^1(\mathcal S^n)$--function~$f$ can be expanded into Laplace series of hyperspherical harmonics by
$$
S(f;x)\sim\sum_{l=0}^\infty Y_l(f;x),
$$
where $Y_l(f;x)$ is given by
\begin{equation*}\begin{split}
Y_l(f;x)&=\frac{\Gamma(\lambda)(\lambda+l)}{2\pi^{\lambda+1}}\int_{\mathcal S^n}f(y)\,C_l^\lambda(x\cdot y)\,d\sigma(y)\\
&=\frac{\lambda+l}{\lambda\Sigma_n}\int_{\mathcal S^n}f(y)\,C_l^\lambda(x\cdot y)\,d\sigma(y).
\end{split}\end{equation*}
For zonal functions we obtain by~\eqref{eq:FH} the representation
$$
Y_l(f;t)=\widehat f(l)\,C_l^\lambda(t),\qquad t=\cos\theta,
$$
with Gegenbauer coefficients
$$
\widehat f(l)=c(l,\lambda)\int_{-1}^1 f(t)\,C_l^\lambda(t)\left(1-t^2\right)^{\lambda-1/2}dt,
$$
where
\begin{equation*}
c(l,\lambda)=\frac{2^{2\lambda-1}\Gamma(l+1)(\lambda+l)\Gamma^2(\lambda)}{\pi\,\Gamma(2\lambda+l)}
   =\frac{\Gamma(l+1)(\lambda+l)\Gamma(\lambda)\Gamma(2\lambda)}{\sqrt\pi\,\Gamma(2\lambda+l)\Gamma(\lambda+1/2)},
\end{equation*}
compare \cite[p. 207]{LV}.
The series
\begin{equation}\label{eq:Gegenbauer_expansion}
\sum_{l=0}^\infty\widehat f(l)\,C_l^\lambda(t)
\end{equation}
is called Gegenbauer expansion of~$f$.

For $f,g\in\mathcal L^1(\mathcal S^n)$, $g$ zonal, their convolution $f\ast g$ is defined by
\begin{equation*}
(f\ast g)(x)=\frac{1}{\Sigma_n}\int_{\mathcal S^n}f(y)\,g(x\cdot y)\,d\sigma(y).
\end{equation*}
With this notation we have
$$
Y_l(f;x)=\frac{\lambda+l}{\lambda}\,\bigl(f\ast C_l^\lambda\bigr)(x),
$$
hence, the function $\mathcal K_l=\frac{\lambda+l}{\lambda}\,C_l^\lambda$ is the reproducing kernel for~$\mathcal H_l(\mathcal S^n)$, and Funk--Hecke formula can be written as
$$
Y_l\ast f=\frac{\lambda}{\lambda+l}\,\widehat f(l)\,Y_l.
$$

Further, any function $f\in\mathcal L^2(\mathcal S^n)$ has a unique representation as a mean--convergent series
\begin{equation}\label{eq:Fourier_series}
f(x)=\sum_{l=0}^\infty\sum_{k\in\mathcal M_{n-1}(l)} a_l^k\,Y_l^k(x),\qquad x\in\mathcal S^n,
\end{equation}
where
$$
a_l^k=a_l^k(f)=\frac{1}{\Sigma_n}\int_{\mathcal S^n}\overline{Y_l^k(x)}\,f(x)\,d\sigma(x)=\left<Y_l^k,f\right>,
$$
for proof cf.~\cite{Vilenkin}. In analogy to the two-dimensional case, we call $a_l^k$ the Fourier coefficients of the function~$f$. Convolution with a zonal function can be then written as
$$
f\ast g=\sum_{l=0}^\infty\sum_{k\in\mathcal M_{n-1}(l)} \frac{\lambda}{\lambda+l}\,a_l^k(f)\,\widehat g(l)\,Y_l^k
$$
and for zonal functions the following relation
$$
\widehat f(l)=A_l^0\cdot a_l^0(f)
$$
between Fourier and Gegenbauer coefficients holds with~$A_l^0$ given by~\eqref{eq:Al0} (compare formulae \eqref{eq:spherical_harmonics}, \eqref{eq:Gegenbauer_expansion}, and~\eqref{eq:Fourier_series}).

The set of rotations of~$\mathbb R^{n+1}$ is denoted by $SO(n+1)$. It is isomorphic to the set of square matrices of degree $n+1$ with determinant~$1$. The $n$--dimensional sphere can be identified with the class of left cosets of $SO(n+1)$ mod $SO(n)$,
$$
\mathcal S^n=SO(n+1)/SO(n),
$$
and $T$ is the regular representation of~$SO(n+1)$ in~$\mathcal L^2(\mathcal S^n)$:
$$
T(\Upsilon)f(x)=f(\Upsilon^{-1}x),
$$
$f\in\mathcal H_l(\mathcal S^n)$, $x\in\mathcal S^n$, $\Upsilon\in SO(n+1)$. The normalized Haar measure on $SO(n+1)$ is denoted by~$d\nu$.

Zonal product of arbitrary $\mathcal L^2(\mathcal S^n)$--functions $f$ and $g$ (introduced in~\cite{EBCK09}) is given by
$$
(f\hat\ast g)(x\cdot y)=\int_{SO(n+1)}f(\Upsilon^{-1}x)\,g(\Upsilon^{-1}y)\,d\nu(\Upsilon),\qquad x,y\in\mathcal S^n,
$$
and it has the representation
\begin{equation}\label{eq:zonal_product_Fc}
(f\hat\ast g)(x\cdot y)=\sum_{l=0}^\infty\sum_{k\in\mathcal M_{n-1}(l)}\frac{a_l^k(f)\,a_l^k(g)}{N(n,l)}\,\frac{\lambda+l}{\lambda}\,C_l^\lambda(x\cdot y).
\end{equation}

For further details on this topic we refer to the textbooks~\cite{Vilenkin} and~\cite{AH12}.

\section{Definition of directional spherical multipole wavelets}\label{sec:definition}

Poisson kernel for the unit sphere located at $\zeta$ inside the unit ball is given by
\begin{equation}\label{eq:Poisson_kernel}
p_\zeta(x)=\frac{1}{\Sigma_n}\frac{1-|\zeta|^2}{|\zeta-x|^{n+1}}=\frac{1}{\Sigma_n}\frac{1-r^2}{(1-2r\cos\theta+r^2)^{(n+1)/2}},
\end{equation}
where
$$
r=|\zeta|<|x|=1
$$
and $\theta$ is the angle between the vectors~$\zeta$ and~$x$, i.e.,
$$
r\cos\theta=\zeta\cdot x.
$$
We assume that $\zeta$ lies on the positive $x_1$--axis, i.e.,
$$
\zeta=r\widehat e.
$$
In this case, $\theta$ is the $\theta_1$--coordinate of~$x$, and Poisson kernel can be written as the series
\begin{equation}\label{eq:PKseries}
p_\zeta(x)=\frac{1}{\Sigma_n}\,\sum_{l=0}^\infty r^l\,\mathcal K_l^\lambda(\cos\theta_1).
\end{equation}
For more details cf.~\cite{SW71} and~\cite{IIN14PW}.

In~\cite{HH09} directional multipole wavelets on a two-dimensional sphere have been introduced. They are defined as derivatives of rotations of Poisson kernel,
$$
g_\rho^{\widehat\varsigma,d}(x)=\rho^d\left.\frac{\partial^d}{\partial\Theta^d}\left(\mathcal R_\zeta(\Theta\widehat\varsigma)\,p_\zeta(x)\right)\right|_{\zeta=e^{-\rho}\widehat e},
$$
where~$\mathcal R_\zeta$ is a rotation acting on~$\zeta$,
$$
\mathcal R_\zeta(\Upsilon)\,p_\zeta(x)=p_{\Upsilon^{-1}\zeta}(x),
$$
with rotation angle~$\Theta$ and rotation axis~$\widehat\varsigma$ parallel to the tangent plane at~$\widehat e$. This concept can be easily applied to the $n$--dimensional case. Without loss of generality we can choose a fixed~$\widehat\varsigma$--axix, since a change of rotation direction in the tangent plane to~$\mathcal S^n$ in~$\widehat e$ results in an $SO(n)$--rotation of the wavelet. For the rest of the paper, we fix~$\widehat\varsigma$ to be the $x_2$--axis.

\begin{bfseries}Notation.\end{bfseries} $\Upsilon_{\!\Theta}$ denotes the rotation of~$\mathbb R^{n+1}$ in the plane $(\widehat e,\widehat\varsigma)$ with rotation angle~$\Theta$. In Cartesian coordinate system, this mapping is given by the matrix
\begin{equation}\label{eq:rotation_matrix}
\left(\begin{array}{ccccc}\cos\Theta&-\sin\Theta&0&\cdots&0\\
\sin\Theta&\cos\Theta&0&\cdots&0\\
0&0&1&\cdots&0\\
\vdots&\vdots&\vdots&\ddots&\vdots\\
0&0&0&\cdots&1\end{array}\right).
\end{equation}

\begin{df}\label{df:direcional_Pwv} Directional Poisson multipole wavelet of order~$d\in\mathbb N$ is defined as
$$
g_\rho^{[d]}(x)=\left.\rho^d\frac{\partial^d}{\partial\Theta^d}\left(\mathcal R_\zeta(\Upsilon_{\!\Theta})\,p_\zeta(x)\right)\right|_{\Theta=0}.
$$
\end{df}

With this notation we have
$$
\mathcal R_\zeta(\Upsilon_{\!\Theta})\,p_\zeta(x)=p_{\Upsilon_{\!\Theta}^{-1}\zeta}(x),
$$
and since $|\Upsilon_{\!\Theta}^{-1}\zeta|=|\zeta|$ and $|\Upsilon_{\!\Theta}^{-1}\zeta-x|=|\zeta-\Upsilon_{\!\Theta} x|$, the rotated kernel is equal to
$$
\mathcal R_\zeta(\Upsilon_{\!\Theta})\,p_\zeta(x)=p_{\zeta}(\Upsilon_{\!\Theta} x)=T(\Upsilon_{\!\Theta}^{-1})\,p_\zeta(x),
$$
cf. formula~\eqref{eq:Poisson_kernel}.

\section{Derivatives of rotations}\label{sec:derivatives_of_rotations}

For further investigation of  Poisson wavelets we need to compute the derivatives of rotation of a single zonal hyperspherical harmonic of degree~$l$ for $\Theta=0$. First we show that thay can be computed recursively, as derivatives of each other, i.e., the evaluation for $\Theta=0$ can be done  after each differentiation.

\begin{lem}\label{lem:recursive_derivation}Let~$Y$ be a zonal function and denote by~$Y^{(d)}$ the derivative of its rotation,
$$
Y^{(d)}(\Theta,x)=\frac{\partial^d}{\partial\Theta^d}Y(\Upsilon_{\!\Theta}x).
$$
Then,
$$
Y^{(d+1)}(0,x)=\frac{\partial}{\partial\Theta} Y^{(d)}(0,\Upsilon_{\!\Theta}x).
$$
\end{lem}

\begin{bfseries}Proof.\end{bfseries}
Functions~$Y^{(d)}$ satisfy
\begin{align*}
Y^{(d)}(\Theta,\Upsilon_{\!\Xi}x)&=\frac{\partial^d}{\partial\Theta^d}Y(\Upsilon_{(\Theta+\Xi)}x)
   =\frac{\partial^d}{\partial\Theta^d}\left[T(\Upsilon_\Xi^{-1})\,Y(\Upsilon_{\!\Theta}x)\right]\\
&=T(\Upsilon_\Xi^{-1})\,\frac{\partial^d}{\partial\Theta^d}Y(\Upsilon_{\!\Theta}x)=T(\Upsilon_\Xi^{-1})\,Y^{(d)}(\Theta,x).
\end{align*}
By change of variables $\Theta=\Theta_1+\Theta_2$, $\Theta_1$ constant, we obtain
\begin{align*}
Y^{(d)}&(\Theta,x)=\partial_{\Theta_2}^d\,Y(\Upsilon_{\!(\Theta_1+\Theta_2)}x)\\
&=\partial_{\Theta_2}^d\,\left[T(\Upsilon_{\!\Theta_1}^{-1})\,Y(\Upsilon_{\Theta_2}x)\right]=T(\Upsilon_{\!\Theta_1}^{-1})\,Y^{(d)}(\Theta_2,x)
\end{align*}
This relation is valid for any pair of variables~$\Theta_1$ and $\Theta_2$ which sum is equal to~$\Theta$. In particular,
$$
Y^{(d)}(\Theta,x)=T(\Upsilon_{\!\Theta}^{-1})\,Y^{(d)}(0,x)=Y^{(d)}(0,\Upsilon_{\!\Theta}x),
$$
and, hence,
$$
Y^{(d+1)}(0,x)=\frac{\partial}{\partial\Theta}\,Y^{(d)}(0,\Upsilon_{\!\Theta}x)\bigr|_{\Theta=0}.
$$
\hfill$\Box$

Using Lemma~\ref{lem:recursive_derivation} we can find Fourier series representations of derivatives of rotations of zonal hyperspherical harmonics.

\begin{lem}\label{lem:derivatives_Ylk}
Let $n\geq3$ and $l\in\mathbb N$ be fixed. Then
\begin{equation}\label{eq:recursion_derTheta}
\left.\frac{\partial}{\partial\Theta}Y_l^{(k_1,0,\dots,0)}(\Upsilon_{\!\Theta}x)\right|_{\Theta=0}
   =\beta_{l,k_1-1}\,Y_l^{(k_1-1,0,\dots,0)}(x)-\beta_{l,k_1}\,Y_l^{(k_1+1,0,\dots,0)}(x)
\end{equation}
for
\begin{equation}\label{eq:beta_l_k1}
\beta_{l,k_1}=\sqrt{\frac{(k_1+1)\,(2\lambda+k_1-1)\,(l-k_1)\,(2\lambda+l+k_1)}{(2\lambda+2k_1-1)\,(2\lambda+2k_1+1)}},
\end{equation}
$k_1=0,1,\dots,l$, and
$$
\beta_{l,-1}=0.
$$
Let $n=2$ and $l\in\mathbb N$ be fixed. Define~$\widetilde{Y}_l^k$ as
$$
\widetilde{Y}_l^k=\begin{cases}Y_l^{-k}+Y_l^k&\text{for }k=0,1,\dots,l,\\0&\text{otherwise.}\end{cases}
$$
Then
\begin{equation}\label{eq:recursion_derThetaY0n2_ver1}
\left.\frac{\partial}{\partial\Theta}\,\widetilde{Y}_l^0(\Upsilon_{\!\Theta}x)\right|_{\Theta=0}=-\sqrt{l(l+1)}\,\widetilde{Y}_l^1(x)
\end{equation}
and
\begin{equation}\label{eq:recursion_derThetan2}
\left.\frac{\partial}{\partial\Theta}\,\widetilde{Y}_l^k(\Upsilon_{\!\Theta}x)\right|_{\Theta=0}=\beta_{l,k-1}\,\widetilde{Y}_l^{k-1}(x)-\beta_{l,k}\,\widetilde{Y}_l^{k+1}(x)
\end{equation}
for $k=1,2,\dots,l$.
\end{lem}

\begin{bfseries}Remark.\end{bfseries} For $\lambda=\frac{1}{2}$ \eqref{eq:beta_l_k1} simplifies to
$$
\beta_{l,k}=\frac{\sqrt{(l-k)\,(l+k+1)}}{2}.
$$
With this formula, equality \eqref{eq:recursion_derThetaY0n2_ver1} becomes
\begin{equation}\label{eq:recursion_derThetaY0n2}
\left.\frac{\partial}{\partial\Theta}\,\widetilde{Y}_l^0(\Upsilon_{\!\Theta}x)\right|_{\Theta=0}=-2\beta_{l,0}\,\widetilde{Y}_l^1(x).
\end{equation}

\begin{bfseries}Proof.\end{bfseries} Denote by $\widetilde\theta_1$, $\widetilde\theta_2$, ..., $\widetilde\varphi$ the spherical coordinates of~$\Upsilon_{\!\Theta}x$. By~\eqref{eq:spherical_coordinates} and~\eqref{eq:rotation_matrix} the following relation holds:
\begin{equation*}\begin{split}
\left(\begin{array}{c}\cos\widetilde\theta_1\\\sin\widetilde\theta_1\cos\widetilde\theta_2\end{array}\right)
   &=\left(\begin{array}{cc}\cos\Theta&-\sin\Theta\\\sin\Theta&\cos\Theta\end{array}\right)\left(\begin{array}{c}\cos\theta_1\\\sin\theta_1\cos\theta_2\end{array}\right)\\
&=\left(\begin{array}{c}\cos\Theta\cos\theta_1-\sin\Theta\sin\theta_1\cos\theta_2\\\sin\Theta\cos\theta_1+\cos\Theta\sin\theta_1\cos\theta_2\end{array}\right).
\end{split}\end{equation*}
It can be seen from this representation that
\begin{equation}\label{eq:derivative_cos1}
\frac{\partial}{\partial\Theta}\cos\widetilde\theta_1=-\sin\widetilde\theta_1\cos\widetilde\theta_2
\end{equation}
and
\begin{equation}\label{eq:derivative_sin1cos2}
\frac{\partial}{\partial\Theta}\,(\sin\widetilde\theta_1\cos\widetilde\theta_2)=\cos\widetilde\theta_1.
\end{equation}
It follows from~\eqref{eq:derivative_cos1} and
$$
\sin\widetilde\theta_1\frac{\partial}{\partial\Theta}\sin\widetilde\theta_1+\cos\widetilde\theta_1\frac{\partial}{\partial\Theta}\cos\widetilde\theta_1=0
$$
that
\begin{equation}\label{eq:derivative_sin1}
\frac{\partial}{\partial\Theta}\sin\widetilde\theta_1=\cos\widetilde\theta_1\cos\widetilde\theta_2
\end{equation}
for $\widetilde\theta_1\in(0,\pi)$, and by continuity of the derivative, for all $\widetilde\theta_1\in[0,\pi]$. In a similar way we obtain from~\eqref{eq:derivative_sin1cos2} and~\eqref{eq:derivative_sin1}
\begin{equation}\label{eq:derivative_cos2}
\sin\widetilde\theta_1\frac{\partial}{\partial\Theta}\cos\widetilde\theta_2=\cos\widetilde\theta_1(1-\cos^2\widetilde\theta_2).
\end{equation}

The two-dimensional case must be treated separately, since then exponential functions of~$\varphi$ are involved. The above relations between spherical coordinates of~$x$ and $\Upsilon_{\!\Theta}$ remain valid with the reservation $\theta_2=\varphi$. Further, it follows from the non-invariance of~$x_3$ under rotation~$\Upsilon_{\!\Theta}$ that
\begin{equation}\label{eq:derivative_sin1sin2}
\frac{\partial}{\partial\Theta}\,(\sin\widetilde\theta_1\sin\widetilde\varphi)=0.
\end{equation}
According to~\eqref{eq:spherical_harmonics} and~\eqref{eq:Alk_n2},
\begin{equation}\label{eq:Ylk_n2}
Y_l^{\pm k}(\widetilde\theta_1,\widetilde\varphi)=2^k\Gamma\left(k+\frac{1}{2}\right)\sqrt{\frac{(2l+1)(l-k)!}{\pi\,(l+k)!}}\,
   C_{l-k}^{k+\frac{1}{2}}(\cos\widetilde\theta_1)\sin^k\widetilde\theta_1\,e^{\pm ik\widetilde\varphi}
\end{equation}
for $k=0,1,\dots,l$. If $k\ne0$,
\begin{equation*}\begin{split}
\frac{\partial}{\partial\Theta}&Y_l^{\pm k}(\widetilde\theta_1,\widetilde\varphi)=2^k\Gamma\left(k+\frac{1}{2}\right)\sqrt{\frac{(2l+1)(l-k)!}{\pi\,(l+k)!}}\\
\cdot&\Bigl[\left(C_{l-k}^{k+\frac{1}{2}}\right)^\prime(\cos\widetilde\theta_1)\frac{\partial\cos\widetilde\theta_1}{\partial\Theta}
   \sin^k\widetilde\theta_1(\cos\widetilde\varphi\pm i\sin\widetilde\varphi)^k\\
&+k\,C_{l-k}^{k+\frac{1}{2}}(\cos\widetilde\theta_1)
   \left(\sin\widetilde\theta_1(\cos\widetilde\varphi\pm i\sin\widetilde\varphi)\right)^{k-1}
   \frac{\partial(\sin\widetilde\theta_1\cos\widetilde\varphi\pm i\sin\widetilde\varphi)}{\partial\Theta}\Bigr]
\end{split}\end{equation*}
Substitute expressions~\eqref{eq:derivative_cos1}, \eqref{eq:derivative_sin1cos2} and~\eqref{eq:derivative_sin1sin2} for the derivatives of trigonometric functions and in the case $k\ne\pm l$ use formula~\eqref{eq:derGg} for the derivative of Gegenbauer polynomial and the fact that $\Upsilon_{\!\Theta}x=x$ for $\Theta=0$ to obtain
\begin{equation*}\begin{split}
&\sqrt{\frac{\pi\,(l+k)!}{(2l+1)(l-k)!}}\left.\frac{1}{2^k\Gamma(k+\frac{1}{2})}\,\frac{\partial Y_l^{\pm k}(\widetilde\theta_1,\widetilde\varphi)}{\partial\Theta}\right|_{\Theta=0}\\
&\qquad=-(2k+1)\,C_{l-k-1}^{k+\frac{3}{2}}(\cos\theta_1)\sin^{k+1}\theta_1\cos\varphi\,(\cos\varphi\pm i\sin\varphi)^k\\
&\qquad+k\,C_{l-k}^{k+\frac{1}{2}}(\cos\theta_1)\sin^{k-1}\theta_1(\cos\varphi\pm i\sin\varphi)^{k-1}\cos\theta_1,
\end{split}\end{equation*}
It is straightforward to verify that
$$
2\cos\varphi\,(\cos\varphi\pm i\sin\varphi)=(\cos\varphi\pm i\sin\varphi)^2+1.
$$
Thus,
\begin{equation}\label{eq:derivative_Ylk}\begin{split}
&\sqrt{\frac{\pi\,(l+k)!}{(2l+1)(l-k)!}}\left.\frac{1}{2^k\Gamma(k+\frac{1}{2})}\,\frac{\partial Y_l^{\pm k}(\widetilde\theta_1,\widetilde\varphi)}{\partial\Theta}\right|_{\Theta=0}\\
&\qquad=-\left(k+\frac{1}{2}\right)C_{l-k-1}^{k+\frac{3}{2}}(\cos\theta_1)\sin^{k+1}\theta_1(\cos\varphi\pm i\sin\varphi)^{k+1}\\
&\qquad+\left[k\cos\theta_1C_{l-k}^{k+\frac{1}{2}}(\cos\theta_1)-\left(k+\frac{1}{2}\right)\sin^2\theta_1C_{l-k-1}^{k+\frac{3}{2}}(\cos\theta_1)\right]\\
&\qquad\quad\cdot\sin^{k-1}\theta_1(\cos\varphi\pm i\sin\varphi)^{k-1}.
\end{split}\end{equation}
The expression in brackets is a polynomial in~$\cos\theta_1$ of degree \mbox{$l-k+1$}. Since the spaces~$\mathcal H_l$, $l\in\mathbb N_0$, are closed, the derivative $\frac{\partial}{\partial\Theta}Y_l(\Upsilon_{\!\Theta}x)$, $Y_l\in\mathcal H_l$, belongs to~$\mathcal H_l$, this polynomial is a multiplicity of $C_{l-(k-1)}^{k+\frac{1}{2}-1}(\cos\theta_1)$ (otherwise the second term on the right-hand-side of~\eqref{eq:derivative_Ylk} would contain a nonzero sum of functions belonging to~$\mathcal H_{l^\prime}$ for $l^\prime< l$). Using representation~\eqref{eq:C_jako_szereg} and comparing coefficients at~$\cos^{l-k+1}\theta_1$ in the regarded polynomial from~\eqref{eq:derivative_Ylk} and~$C_{l-(k-1)}^{k+\frac{1}{2}-1}(\cos\theta_1)$, we find out that the proportionality constant is equal to
$$
\frac{(l-k+1)(l+k)}{2(2k-1)}.
$$
Now, the expressions on the right-hand side of~\eqref{eq:derivative_Ylk} may be replaced by hyperspherical harmonics~$Y_l^{\pm(k+1)}$ and~$Y_l^{\pm(k-1)}$ given by~\eqref{eq:Ylk_n2} and we obtain
\begin{equation*}\begin{split}
\left.\frac{\partial Y_l^{\pm k}(\widetilde\theta_1,\widetilde\varphi)}{\partial\Theta}\right|_{\Theta=0}&=-\frac{\sqrt{(l-k)(l+k+1)}}{2}\,Y_l^{\pm(k+1)}(\theta_1,\varphi)\\
&+\frac{\sqrt{(l-k+1)(l+k)}}{2}\,Y_l^{\pm(k-1)}(\theta_1,\varphi).
\end{split}\end{equation*}

Analogous considerations yield
$$
\left.\frac{\partial}{\partial\Theta}\,Y_l^0(\Upsilon_{\!\Theta}x)\right|_{\Theta=0}=-\frac{\sqrt{l(l+1)}}{2}\left[Y_l^1(x)+Y_l^{-1}(x)\right]
$$
and
$$
\left.\frac{\partial}{\partial\Theta}\,Y_l^{\pm l}(\Upsilon_{\!\Theta}x)\right|_{\Theta=0}=\sqrt{\frac{l}{2}}\,Y_l^{\pm(l-1)}(x).
$$

For $n\geq3$ we have by~\eqref{eq:spherical_harmonics}
$$
Y_l^{(k_1,0,\dots,0)}(\Upsilon_{\!\Theta}x)
   =A_l^{(k_1,0,\dots,0)}\,C_{l-k_1}^{\lambda+k_1}(\cos\widetilde\theta_1)\sin^{k_1}\widetilde\theta_1\,C_{k_1}^{\lambda-\frac{1}{2}}(\cos\widetilde\theta_2).
$$

Consequently,  for $k_1\ne0$,
\begin{align*}
\frac{\partial}{\partial\Theta}&\frac{Y_l^{(k_1,0,\dots,0)}(\Upsilon_{\!\Theta}x)}{A_l^{(k_1,0,\dots,0)}}\\
&=(C_{l-k_1}^{\lambda+k_1})^\prime(\cos\widetilde\theta_1)\frac{\cos\widetilde\theta_1}{\partial\Theta}
   \sin^{k_1}\widetilde\theta_1\,C_{k_1}^{\lambda-\frac{1}{2}}(\cos\widetilde\theta_2)\\
&+k_1\,C_{l-k_1}^{\lambda+k_1}(\cos\widetilde\theta_1)\sin^{k_1-1}\widetilde\theta_1\frac{\partial\sin\widetilde\theta_1}{\partial\Theta}
   \,C_{k_1}^{\lambda-\frac{1}{2}}(\cos\widetilde\theta_2)\\
&+C_{l-k_1}^{\lambda+k_1}(\cos\widetilde\theta_1)\sin^{k_1}\widetilde\theta_1\,(C_{k_1}^{\lambda-\frac{1}{2}})^\prime(\cos\widetilde\theta_2)
   \frac{\partial\cos\widetilde\theta_2}{\partial\Theta}.
\end{align*}
We substitute expressions~\eqref{eq:derivative_cos1}, \eqref{eq:derivative_sin1} and~\eqref{eq:derivative_cos2} for the derivatives of trigonometric functions, use the fact that $\Upsilon_{\!\Theta}x=x$ for $\Theta=0$, and obtain by~\eqref{eq:derGg}
\begin{equation}\label{eq:derivativeYlk_Theta0}\begin{split}
\frac{\partial}{\partial\Theta}&\left.\frac{Y_l^{(k_1,0,\dots,0)}(\Upsilon_{\!\Theta}x)}{A_l^{(k_1,0,\dots,0)}}\right|_{\Theta=0}\\
&=-2(\lambda+k_1)\,C_{l-k_1-1}^{\lambda+k_1+1}(\cos\theta_1)\sin^{k_1+1}\theta_1\cos\theta_2\,C_{k_1}^{\lambda-\frac{1}{2}}(\cos\theta_2)\\
&+k_1\,C_{l-k_1}^{\lambda+k_1}(\cos\theta_1)\sin^{k_1-1}\theta_1\cos\theta_1\cos\theta_2\,C_{k_1}^{\lambda-\frac{1}{2}}(\cos\theta_2)\\
&+(2\lambda-1)\,C_{l-k_1}^{\lambda+k_1}(\cos\theta_1)\sin^{k_1-1}\theta_1\,C_{k_1-1}^{\lambda+\frac{1}{2}}(\cos\theta_2)\cos\theta_1(1-\cos^2\theta_2)
\end{split}\end{equation}
for $k_1\ne l$. Formula~\eqref{eq:recGg1} yields the following representation for the first term in~\eqref{eq:derivativeYlk_Theta0}:
\begin{equation*}\begin{split}
-&\frac{2(\lambda+k_1)}{2\lambda+2k_1-1}\,C_{l-k_1-1}^{\lambda+k_1+1}(\cos\theta_1)\sin^{k_1+1}\theta_1\\
&\cdot\left[(k_1+1)\,C_{k_1+1}^{\lambda-\frac{1}{2}}(\cos\theta_2)+(2\lambda+k_1-2)\,C_{k_1-1}^{\lambda-\frac{1}{2}}(\cos\theta_2)\right]\\
&=-\frac{2(\lambda+k_1)(k_1+1)}{2\lambda+2k_1-1}\,C_{l-k_1-1}^{\lambda+k_1+1}(\cos\theta_1)\sin^{k_1+1}\theta_1\,C_{k_1+1}^{\lambda-\frac{1}{2}}(\cos\theta_2)\\
&+\frac{2(\lambda+k_1)(2\lambda+k_1-2)}{2\lambda+2k_1-1}\,(\cos^2\theta_1-1)\,
   C_{l-k_1-1}^{\lambda+k_1+1}(\cos\theta_1)\sin^{k_1-1}\theta_1\,C_{k_1-1}^{\lambda-\frac{1}{2}}(\cos\theta_2)
\end{split}\end{equation*}

The common factor in the second and the third summand on the right-hand-side~\eqref{eq:derivativeYlk_Theta0} can be written as
\begin{equation*}\begin{split}
\cos\theta_1&\,C_{l-k_1}^{\lambda+k_1}(\cos\theta_1)\sin^{k_1-1}\theta_1\\
&=\left[\frac{l-k_1+1}{2(\lambda+k_1-1)}\,C_{l-k_1+1}^{\lambda+k_1-1}(\cos\theta_1)+C_{l-k_1-1}^{\lambda+k_1}(\cos\theta_1)\right]\sin^{k_1-1}\theta_1
\end{split}\end{equation*}
according to~\eqref{eq:recGg2}. The remaining terms in the second and the third summand in~\eqref{eq:derivativeYlk_Theta0},
$$
k_1\,\cos\theta_2\,C_{k_1}^{\lambda-\frac{1}{2}}(\cos\theta_2)+(2\lambda-1)(1-\cos^2\theta_2)\,C_{k_1-1}^{\lambda+\frac{1}{2}}(\cos\theta_2),
$$
can be with formula~\eqref{eq:recGg3} rewritten as
\begin{equation*}\begin{split}
k_1\,\cos\theta_2\,&C_{k_1}^{\lambda-\frac{1}{2}}(\cos\theta_2)
   +\left[(2\lambda+k_1-1)\,\cos\theta_2\,C_{k_1}^{\lambda-\frac{1}{2}}(\cos\theta_2)-(k_1+1)\,C_{k_1+1}^{\lambda-\frac{1}{2}}(\cos\theta_2)\right]\\
&=(2\lambda+2k_1-1)\,\cos\theta_2\,C_{k_1}^{\lambda-\frac{1}{2}}(\cos\theta_2)-(k_1+1)\,C_{k_1+1}^{\lambda-\frac{1}{2}}(\cos\theta_2),
\end{split}\end{equation*}
and further, according to formula~\eqref{eq:recGg1}, they are equal to
$$
(2\lambda+k_1-2)\,C_{k_1-1}^{\lambda-\frac{1}{2}}(\cos\theta_2).
$$
Consequently, from \eqref{eq:derivativeYlk_Theta0} we obtain
\begin{equation*}\begin{split}
\frac{\partial}{\partial\Theta}&\left.\frac{Y_l^{(k_1,0,\dots,0)}(\Upsilon_{\!\Theta}x)}{A_l^{(k_1,0,\dots,0)}}\right|_{\Theta=0}\\
&=-\frac{2(\lambda+k_1)(k_1+1)}{2\lambda+2k_1-1}\,C_{l-k_1-1}^{\lambda+k_1+1}(\cos\theta_1)\sin^{k_1+1}\theta_1\,C_{k_1+1}^{\lambda-\frac{1}{2}}(\cos\theta_2)\\
&+\left[\frac{2(\lambda+k_1)}{2\lambda+2k_1-1}\,(\cos^2\theta_1-1)\,C_{l-k_1-1}^{\lambda+k_1+1}(\cos\theta_1)\right.\\
&\quad+\left.\frac{(l-k_1+1)}{2(\lambda+k_1-1)}\,C_{l-k_1+1}^{\lambda+k_1-1}(\cos\theta_1)+C_{l-k_1-1}^{\lambda+k_1}(\cos\theta_1)\right]\sin^{k_1-1}\theta_1\\
&\quad\cdot(2\lambda+k_1-2)\,C_{k_1-1}^{\lambda-\frac{1}{2}}(\cos\theta_2).
\end{split}\end{equation*}
The expression in brackets is a polynomial in~$\cos\theta_1$ of degree \mbox{$l-k+1$}.  Since $\frac{\partial}{\partial\Theta}Y_l(\Upsilon_{\!\Theta}x)$ belongs to~$\mathcal H_l$ for $Y_l\in\mathcal H_l$, this polynomial is a multiplicity of~$C_{l-k_1+1}^{\lambda+k_1-1}$. A comparison of $l-k_1+1$--st coefficients yields the proportionality constant to be equal to
\begin{equation*}
\frac{(l-k_1+1)(2\lambda+l+k_1-1)}{2(\lambda+k_1-1)(2\lambda+2k_1-1)}.
\end{equation*}
Consequently,
\begin{align}
\frac{\partial}{\partial\Theta}&Y_l^{(k_1,0,\dots,0)}(\Upsilon_{\!\Theta}x)
   =-\frac{2(\lambda+k_1)\,(k_1+1)}{2\lambda+2k_1-1}\frac{A_l^{(k_1,0,\dots,0)}}{A_l^{(k_1+1,0,\dots,0)}}\,Y_l^{(k_1+1,0,\dots,0)}(x)\label{eq:beta_k}\\
&+\frac{(l-k_1+1)\,(2\lambda+l+k_1-1)\,(2\lambda+k_1-2)}{2(\lambda+k_1-1)\,(2\lambda+2k_1-1)}\frac{A_l^{(k_1,0,\dots,0)}}{A_l^{(k_1-1,0,\dots,0)}}\,
   Y_l^{(k_1-1,0,\dots,0)}(x)\notag\\
&=-\sqrt{\frac{(k_1+1)\,(2\lambda+k_1-1)\,(l-k_1)\,(2\lambda+l+k_1)}{(2\lambda+2k_1-1)\,(2\lambda+2k_1+1)}}\,Y_l^{(k_1+1,0,\dots,0)}(x)\notag\\
&+\sqrt{\frac{k_1(2\lambda+k_1-2)\,(l-k_1+1)\,(2\lambda+l+k_1-1)}{(2\lambda+2k_1-3)\,(2\lambda+2k_1-1)}}\,Y_l^{(k_1-1,0,\dots,0)}(x).\notag
\end{align}
In an analogous way we obtain
$$
\left.\frac{\partial}{\partial\Theta}Y_l^0(\Upsilon_{\!\Theta}x)\right|_{\Theta=0}=-\sqrt{\frac{(2\lambda+l)\,l}{2\lambda+1}}\,Y_l^1(x),
$$
and
$$
\left.\frac{\partial}{\partial\Theta}Y_l^l(\Upsilon_{\!\Theta}x)\right|_{\Theta=0}=\sqrt{\frac{l\,(2\lambda+l-2)}{2\lambda+2l-3}}\,Y_l^{l-1}(x).
$$
\hfill$\Box$

As a consequence, from formulae~\eqref{eq:recursion_derTheta}, \eqref{eq:recursion_derThetaY0n2} and \eqref{eq:recursion_derThetan2} we obtain the following representation of derivatives of functions.

\begin{thm}\label{thm:series_representation}Let a zonal function
$$
f=\sum_{l=0}^\infty a_l^0(f)\,Y_l^0
$$
be given. Then
\begin{equation}\label{eq:fd_series_n3}
f^{(d)}:=\left.\frac{\partial^{d}}{\partial\Theta^{d}}\,\left(f(\Upsilon_\Theta x)\right)\right|_{\Theta=0}
   =\sum_{l=0}^\infty\sum_{j=0}^{\left[\frac{d}{2}\right]}a_l^{2j+d_{\text{mod}2}}(f^{(d)})Y_l^{(2j+d_{\text{mod}2},0,\dots,0)}(x)
\end{equation}
for $n\geq3$ or
\begin{equation}\label{eq:fd_series_n2}
f^{(d)}=\sum_{l=0}^\infty\sum_{j=0}^{\left[\frac{d}{2}\right]}a_l^{2j+d_{\text{mod}2}}(f^{(d)})\left(Y_l^{2j+d_{\text{mod}2}}+Y_l^{-(2j+d_{\text{mod}2})}\right)
\end{equation}
for $n=2$ with coefficients~$a_l^k(f^{(d)})$, obtained recursively via
\begin{align*}
a_l^{2j+1}(f^{(d)})&=\beta_{l,2j+1}a_l^{2j+2}(f^{(d-1)})-\beta_{l,2j}a_l^{2j}(f^{(d-1)}),\\
a_l^{2j}(f^{(d)})&=0
\end{align*}
for an odd~$d$ and
\begin{align*}
a_l^0(f^{(d)})&=-\beta_{l,0}\,a_l^1(f^{(d-1)})\qquad(\text{for }n\geq3),\\
a_l^0(f^{(d)})&=-2\beta_{l,0}\,a_l^1(f^{(d-1)})\qquad(\text{for }n=2),\\
a_l^{2j}(f^{(d)})&=\beta_{l,2j-1}\,a_l^{2j-1}(f^{(d-1)})-\beta_{l,2j}\,a_l^{2j+1}(f^{(d-1)}),\\
a_l^{2j+1}(f^{(d)})&=0
\end{align*}
for an even~$d$, where $\beta_{l,k}$ are defined as in Lemma~\ref{lem:derivatives_Ylk}.
\end{thm}

\begin{bfseries}Remark~1. \end{bfseries}We adapt the convention $Y_l^{(k_1,0,\dots,0)}=0$ for $k_1\geq l$. This can be also interpreted as $a_l^{k_1}(f^{(d)})=0$ for $k_1\geq l$.

\begin{bfseries}Remark~2.\end{bfseries} In~\cite{HH09} this recursion is represented by a matrix multiplication. Unfortunately, it is disregarded that the infinitesimal rotation operator acts on~$Y_l^0$ in a slightly different way than creation and anihilation, cf.~\eqref{eq:recursion_derThetaY0n2}. However, only a constant must be corrected in the expressions representing~$\widehat g_\rho^n(l,m)$, $E_l^{[n]}(a)$, and $c_g(l)$.

\begin{cor}\label{cor:f(1)explicit}Let a zonal function
$$
f(x)=\sum_{l=0}^\infty \widehat f(l)\,C_l^\lambda(\cos\theta_1)
$$
for $x=(\theta_1,\theta_2,\dots,\theta_{n-1},\varphi)$ in spherical coordinates be given. Then
\begin{equation}\label{eq:f1}
f^{(1)}(x):=\left.\frac{\partial}{\partial\Theta}f(\Upsilon_\Theta x)\right|_{\Theta=0}=\frac{\partial}{\partial\theta_1}f(x)\cdot\cos\theta_2.
\end{equation}
\end{cor}

\begin{bfseries}Proof.\end{bfseries} For $n\geq3$ from Theorem~\ref{thm:series_representation} and \eqref{eq:spherical_harmonics} we have
\begin{equation*}
f^{(1)}=-\sum_{l=1}^\infty\beta_{l,0}\,\frac{\widehat f(l)}{A_l^0}\,Y_l^1.
\end{equation*}
Use~\eqref{eq:spherical_harmonics} for~$Y_l^1$ and substitute $k_1=0$ to the factor occurring in~\eqref{eq:beta_k} to express~$\beta_{l,0}$ to obtain
\begin{equation*}\begin{split}
f^{(1)}(x)&=-\sum_{l=1}^\infty\frac{2\lambda}{2\lambda-1}\frac{A_l^0}{A_l^{(1,0,\dots,0)}}\,\frac{\widehat f(l)}{A_l^0}
   \cdot A_l^{(1,0,\dots,0)}\,C_{l-1}^{\lambda+1}(\cos\theta_1)\,\sin\theta_1\,C_1^{\lambda-\frac{1}{2}}(\cos\theta_2)\\
&=\frac{1}{2\lambda-1}\sum_{l=1}^\infty \widehat f(l)\,\frac{dC_l^\lambda(\cos\theta_1)}{d\theta_1}\,\,C_1^{\lambda-\frac{1}{2}}(\cos\theta_2),
\end{split}\end{equation*}
according to~\eqref{eq:derGg}. \eqref{eq:f1} follows by~\cite[formula~8.930.2]{GR} (compare also~\eqref{eq:C_jako_szereg}). For $n=2$ the proof is analogous.\hfill$\Box$

\begin{exa}\label{exa:g_rho_1}An explicit representation of~$g_\rho^{[1]}$ is given by
\begin{equation}\label{eq:g_rho_1_explicit}
g_\rho^{[1]}(x)=-\frac{2\rho\,(\lambda+1)\,r\,(1-r^2)\,\sin\theta_1}{\Sigma_n(1-2r\cos\theta_1+r^2)^{\lambda+2}}\cdot\cos\theta_2,\qquad r=e^{-\rho}.
\end{equation}
\end{exa}

\begin{cor}\label{cor:structure_alpha}Coefficients~$a_l^j(f^{(d)})$ in~\eqref{eq:fd_series_n3} and~\eqref{eq:fd_series_n2} have the following structure:
\begin{equation}\label{eq:structure_aljf}
a_l^j(f^{(d)})=\left(\prod_{\iota=0}^{j-1}\beta_{l,\iota}\right)\cdot p_{d,j}\,(\beta_{l,0}^2,\beta_{l,1}^2,\dots,\beta_{l,\frac{d-j}{2}}^2)\cdot a_l^j(f)
\end{equation}
for $d$, $j$ with the same parity, where~$p_{d,j}$ is a polynomial of degree~$\frac{d-j}{2}$, and
\begin{equation}\label{eq:aljf_equal_zero}
a_l^j(f^{(d)})=0
\end{equation}
otherwise.
\end{cor}

\begin{bfseries}Proof. \end{bfseries}The statement holds for $d=0$ and $d=1$. Suppose, it holds for some even~$d$. Then, for $d+1$ we have
\begin{align*}
a_l^{2j+1}(f^{(d+1)})&=\beta_{l,2j+1}a_l^{2j+2}(f^{(d)})-\beta_{l,2j}a_l^{2j}(f^{(d)})\\
&=\beta_{l,2j+1}\left(\prod_{\iota=0}^{2j+1}\beta_{l,\iota}\right)\cdot p_{d,2j+2}\,(\beta_{l,0}^2,\beta_{l,1}^2,\dots,\beta_{l,\frac{d}{2}-j-1}^2)\cdot a_l^j(f)\\
&-\beta_{l,2j}\left(\prod_{\iota=0}^{2j-1}\beta_{l,\iota}\right)\cdot p_{d,2j}\,(\beta_{l,0}^2,\beta_{l,1}^2,\dots,\beta_{l,\frac{d}{2}-j}^2)\cdot a_l^j(f)\\
&=\left(\prod_{\iota=0}^{2j}\beta_{\iota}\right)\cdot p_{d+1,2j+1}(\beta_{l,0}^2,(\beta_{l,1}^2,\dots,\beta_{l,\frac{d}{2}-j}^2)\cdot a_l^j(f)
\end{align*}
for
$$
p_{d+1,2j+1}(\beta_{0}^2,\dots,\beta_{\frac{d}{2}-j}^2)=\beta_{2j+1}^2\cdot p_{d,2j+2}\,(\beta_{0}^2,\dots,\beta_{\frac{d}{2}-j-1}^2)-p_{d,2j}\,(\beta_{0}^2,\dots,\beta_{\frac{d}{2}-j}^2).
$$
The proof is analogous for odd~$d$'s.\hfill$\Box$

\begin{bfseries}Remark. \end{bfseries}Formula~\eqref{eq:structure_aljf} is valid also for $j>l$, since~$\beta_{l,l}=0$. In this case, $\beta_{l,j}$, $j>l$, has to be understood as an imaginary number.

\begin{prop}\label{prop:well_definiteness}
Directional Poisson wavelets are well--defined $\mathcal L^2$--functions.
\end{prop}

\begin{bfseries}Proof.\end{bfseries} According to~\eqref{eq:Al0},
$$
A_l^0=\mathcal O(l^{\frac{3-n}{2}}),\qquad\text{for }l\to\infty.
$$
It follows from~\eqref{eq:PKseries} and from Corollary~\ref{cor:structure_alpha} that the coefficients of~$g_\rho^{[d]}$ in the representation as a series of hyperspherical harmonics behave for large~$l$ like
$$
l^\alpha r^{-l}
$$
for some rational~$\alpha$'s, thus, they are square summable, and the functions are square integrable.\hfill$\Box$

\section{Wavelets derived from an approximate identity}\label{sec:wavelets}

Directional wavelets on $n$--dimensional spheres derived from an approximate identity were studied in~\cite{EBCK09} and in~\cite{IIN14CWT}. For our purposes we need to modify the definition used in those papers slightly.

\begin{df}\label{def:bilinear_wavelets} Let $\alpha:\mathbb R_+\to\mathbb R_+$ be a weight function. Families $\{\Psi_\rho\}_{\rho\in\mathbb R_+}\subseteq\mathcal L^2(\mathbb S^n)$ and $\{\Omega_\rho\}_{\rho\in\mathbb R_+}\subseteq\mathcal L^2(\mathbb S^n)$ are called an admissible wavelet pair if they satisfy the following conditions:
\begin{enumerate}
\item for $l\in\mathbb{N}_0$
\begin{equation}\label{eq:admbwv1}
\sum_{\kappa=1}^{N(n,l)}\int_0^\infty\overline{a_l^\kappa(\Psi_\rho)}\,a_l^\kappa(\Omega_\rho)\,\alpha(\rho)d\rho=N(n,l),
\end{equation}
\item for $R\in\mathbb{R}_+$ and $x\in\mathbb S^n$
\begin{equation}\label{eq:admbwv2}
\int_{\mathbb S^n}\left|\int_R^\infty(\overline{\Psi_\rho}\hat\ast\Omega_\rho)(x\cdot y)\,\alpha(\rho)d\rho\right|d\sigma(y)\leq\mathfrak c
\end{equation}
with $\mathfrak c$ independent of~$R$.
\end{enumerate}
\end{df}

\begin{df}\label{def:bilinear_wt} Let $\{\Psi_\rho\}_{\rho\in\mathbb R_+}$ and $\{\Omega_\rho\}_{\rho\in\mathbb R_+}$ be an admissible wavelet pair. Then, the spherical wavelet transform
$$
\mathcal W_\Psi\colon\mathcal L^2(\mathcal S^n)\to\mathcal L^2(\mathbb R_+\times SO(n+1))
$$
is defined by
$$
\mathcal W_\Psi f(\rho,\Upsilon)=\frac{1}{\Sigma_n}\int_{\mathcal S^n}\overline{\Psi_\rho(\Upsilon^{-1}x)}\,f(x)\,d\sigma(x).
$$
\end{df}

It can be proven in a similar way as in~\cite{EBCK09} that the wavelet transform is invertible by
$$
f(x)=\int_0^\infty\!\!\int_{SO(n+1)}\mathcal W_\Psi f(\rho,\Upsilon)\,\Omega_\rho(\Upsilon^{-1}x)\,d\nu(\Upsilon)\,\alpha(\rho)d\rho.
$$
Unless $\{\Psi_\rho\}$ and $\{\Omega_\rho\}$ are not equal to each other, the wavelet transform is not an isometry, compare the proof of \cite[Theorem~3.3]{IIN14CWT}.

\begin{bfseries}Remark.\end{bfseries} The wavelet transform can be written as
$$
\mathcal W_\Psi f(\rho,\Upsilon)=\left<\Psi_{\rho,\Upsilon},f\right>=f\ast\overline{\Psi_\rho}(\Upsilon),\qquad\Psi_{\rho,\Upsilon}=\Psi_\rho(\Upsilon^{-1}\circ),\,\Upsilon\in SO(n+1).
$$

Our goal in this section is to prove that for slightly modified Poisson wavelets families of functions exist such that they build admissible wavelet pairs for $\alpha(\rho)=\frac{1}{\rho}$. For~\eqref{eq:admbwv1} to be satisfied the functions need to have the property described in the following Conjecture.

\begin{con}\label{con:linear_combination_of_derivatives}
Let $\mathcal L^2$--functions
$$
f=\sum_{l=0}^\infty a_l^0(f)Y_l^0\qquad\text{and}\qquad g=\sum_{l=0}^\infty a_l^0(g)Y_l^0
$$
be given and denote by~$f^{(d)}$, $g^{(d)}$, $d\in\mathbb N$, derivatives of their rotation,
$$
f^{(d)}(x)=\left.\frac{\partial^d}{\partial\Theta^d}\left(f(\Upsilon_\Theta x)\right)\right|_{\Theta=0},\qquad
g^{(d)}(x)=\left.\frac{\partial^d}{\partial\Theta^d}\left(g(\Upsilon_\Theta x)\right)\right|_{\Theta=0}
$$
$x\in\mathcal S^n$. For any $\mathfrak d$ there exist coefficients $\gamma_0^\mathfrak d$, $\gamma_1^\mathfrak d$, $\dots$, $\gamma_{\mathfrak d-1}^\mathfrak d$, and $\gamma_\mathfrak d^\mathfrak d$, independent of~$f$ and~$g$, such that Fourier coefficients of
$$
F^{\mathfrak d}=\sum_{d=0}^\mathfrak d\gamma_d^\mathfrak d f^{(d)}
   \qquad\text{and}\qquad G^{\mathfrak d}=\sum_{d=0}^\mathfrak d\gamma_d^\mathfrak d g^{(d)}
$$
satisfy
\begin{equation}\label{eq:thesis_conjecture}
\sum_{k_1=0}^l\overline{a_l^{k_1}(F^\mathfrak d)}\,a_l^{k_1}(G^\mathfrak d)=\overline{a_l^0(f)}\,a_l^0(g)\,[l(2\lambda+l)]^\mathfrak d.
\end{equation}
\end{con}

\begin{bfseries}Proof. \end{bfseries}According to Theorem~\ref{thm:series_representation} and Corollary~\ref{cor:structure_alpha} (and the Remark following it),
$$
f^{(d)}=\sum_{l=0}^\infty\sum_{j=0}^\mathfrak d a_l^j(f^{(d)})\,Y_l^j,\qquad g^{(d)}=\sum_{l=0}^\infty\sum_{j=0}^\mathfrak d a_l^j(g^{(d)})\,Y_l^j
$$
with coefficients~$a_l^j$ satisfying~\eqref{eq:structure_aljf}, resp.~\eqref{eq:aljf_equal_zero}. Note that~$\beta_{l,j}^2$ are equal to
\begin{equation}\label{eq:beta_linear}
\beta_{l,j}^2=\frac{(j+1)\,(2\lambda+j-1)}{(2\lambda+2j-1)\,(2\lambda+2j+1)}\,[l(2\lambda+l)-j(2\lambda+j)],
\end{equation}
cf.~\eqref{eq:beta_l_k1}, i.e., they are linear polynomials of the variable~$l(2\lambda+l)$. Projections~$F_l^\mathfrak d$ resp. $G_l^\mathfrak d$ of~$F^{\mathfrak d}$ resp. $G^{\mathfrak d}$ onto~$\mathcal H_l$ are given by
\begin{align*}
F_l^\mathfrak d&=\sum_{d=0}^\mathfrak d\sum_{j=0}^\mathfrak d\gamma_d^\mathfrak d\,a_l^j(f^{(d)})\,Y_l^j
   =\sum_{j=0}^\mathfrak d\sum_{d=0}^\mathfrak d\gamma_d^\mathfrak d\,a_l^j(f^{(d)})\,Y_l^j,\\
G_l^\mathfrak d&=\sum_{j=0}^\mathfrak d\sum_{d=0}^\mathfrak d\gamma_d^\mathfrak d\,a_l^j(g^{(d)})\,Y_l^j.
\end{align*}
For the sum of products of their Fourier coefficients we obtain
\begin{equation}\label{eq:sum_alpha_squared}\begin{split}
\sum_{j=0}^{N(n,l)}\overline{a_l^j(F^\mathfrak d)}\,a_l^j(G^\mathfrak d)
   &=\sum_{j=0}^\mathfrak d\sum_{d=0}^\mathfrak d\sum_{d^\prime=0}^\mathfrak d\overline{\gamma_d^\mathfrak d}\gamma_{d^\prime}^\mathfrak d\overline{a_l^j(f^{(d)})}\,a_l^j(g^{(d^\prime)})\\
   &=\sum_{d=0}^\mathfrak d\sum_{d^\prime=0}^\mathfrak d\overline{\gamma_d^\mathfrak d}\gamma_{d^\prime}^\mathfrak d\sum_{j=0}^\mathfrak d\overline{a_l^j(f^{(d)})}\,a_l^j(g^{(d^\prime)}).
\end{split}\end{equation}
According to~\eqref{eq:structure_aljf}, \eqref{eq:aljf_equal_zero}, and~\eqref{eq:beta_linear}, $\sum_{j=0}^\mathfrak d\overline{a_l^j(f^{(d)})}\,a_l^j(g^{(d^\prime)})$ is a polynomial in \mbox{$l(2\lambda+l)$} of degree $j+\frac{d-j}{2}+\frac{d^\prime-j}{2}=\frac{d+d^\prime}{2}$ multiplied by~$\overline{a_l^0(f)}\,a_l^0(g)$ if~$d$ and~$d^\prime$ are of the same parity, and~$0$ otherwise. Denote it by~$q_{d,d^\prime}$, then
\begin{align*}
\sum_{j=0}^l|\widehat\alpha_{l,j}^\mathfrak d|^2&=|\gamma_\mathfrak d^\mathfrak d|^2q_{\mathfrak d,\mathfrak d}
   +\sum_{\iota=0}^2\overline{\gamma_{\mathfrak d-\iota}^\mathfrak d}\gamma_{\mathfrak d-2+\iota}^\mathfrak d\,q_{\mathfrak d-\iota,\mathfrak d-2+\iota}
   +\sum_{\iota=0}^4\overline{\gamma_{\mathfrak d-\iota}^\mathfrak d}\gamma_{\mathfrak d-4+\iota}^\mathfrak d\,q_{\mathfrak d-\iota,\mathfrak d-2+\iota}\\
&+\dots+|\gamma_0^\mathfrak d|^2\,q_{0,0}.
\end{align*}
The requirement that this expression is equal to~$\overline{a_l^0(f)}\,a_l^0(g)\,[l(2\lambda+l)]^\mathfrak d$ leads to the system of equations
\begin{equation}\label{eq:system_gamma}\begin{split}
\mathfrak c_{\mathfrak d,\mathfrak d}^{(\mathfrak d)}\cdot|\gamma_{\mathfrak d}^\mathfrak d|^2&=\overline{a_l^0(f)}\,a_l^0(g)\\
\mathfrak c_{\mathfrak d,\mathfrak d}^{(\mathfrak d-J)}+\sum_{j=1}^J\sum_{\iota=0}^{2j}
   \mathfrak c_{\mathfrak d-\iota,\mathfrak d-2j+\iota}^{(\mathfrak d-J)}\gamma_{\mathfrak d-\iota}^\mathfrak d\gamma_{\mathfrak d-2j+\iota}^\mathfrak d&=0,
   \qquad J=1,2,\dots,\mathfrak d.
\end{split}\end{equation}
where~$\mathfrak c_{d,d^\prime}^{(J)}$ denote the coefficients of~$q_{d,d^\prime}$,
$$
q_{d,d^\prime}(t)=\sum_{J=0}^{\frac{d+d^\prime}{2}}\mathfrak c_{d,d^\prime}^{(J)}t^J.
$$
The $J^\text{th}$ equation, $J=1,2,\dots,\mathfrak d$, describes vanishing of $(\mathfrak d-J)^\text{th}$ coefficient of polynomial~\eqref{eq:sum_alpha_squared} and $\gamma_j^\mathfrak d=0$ for $j\leq0$. We are not able to compute the Gr\"obner basis of the system~\eqref{eq:system_gamma} in this general setting, hence, we cannot prove that a solution exists.\hfill$\Box$

\begin{exa}
\begin{align*}
F^\mathfrak d=G^\mathfrak d&=f,\\
F^\mathfrak d=G^\mathfrak d&=\sqrt{2\lambda+1}\,f^{(1)},\\
F^\mathfrak d=G^\mathfrak d&=2\,\sqrt{\frac{\lambda(2\lambda+1)}{3}}\,f^{(1)}+\sqrt{\frac{(2\lambda+1)(3+2\lambda)}{3}}f^{(2)},\\
F^\mathfrak d=G^\mathfrak d&=4\,\sqrt{\frac{(-\lambda+1)\lambda(2\lambda+1)^2(3+2\lambda)(5+2\lambda)}{15(2\lambda+1)(3+2\lambda)(5+2\lambda)}}\,f^{(1)}\\
&+2\,\sqrt{\frac{(2\lambda+1)\left[2\sqrt{(-\lambda+1)\lambda(3+2\lambda)(5+2\lambda)}+5\lambda(3+2\lambda)\right]}{15}}\,f^{(2)}\\
&+\sqrt{\frac{(2\lambda+1)(3+2\lambda)(5+2\lambda)}{15}}\,f^{(3)}
\end{align*}
satisfy the thesis of Conjecture~\ref{con:linear_combination_of_derivatives}.
\end{exa}

Now we come to the main result of this section, namely, a construction of admissible wavelet pairs.

\begin{df}\label{df:mdPw}Let $g_\rho^{(0)}=g_\rho$ denote Poisson kernel~$p_\zeta$ located at $\zeta=e^{-\rho}\hat e$ and $h_\rho^{(0)}=h_\rho$ be given by
$$
h_\rho(x)=\frac{1}{\Sigma_n}\sum_{l=0}^\infty e^{-\frac{\rho l^2}{2\lambda}}\,\mathcal K_l^\lambda(\cos\theta_1)
$$
for $x=(\theta_1,\theta_2,\dots,\theta_{n-1},\varphi)$ in spherical coordinates. Further, denote by $g_\rho^{(d)}$ resp. $h_\rho^{(d)}$, $d\in\mathbb N$, directional derivatives of these functions,
$$
g_\rho^{(d)}(x)=\left.\frac{\partial^d}{\partial\Theta^d}\,g_\rho(\Upsilon_\Theta x)\right|_{\Theta=0},\qquad
h_\rho^{(d)}(x)=\left.\frac{\partial^d}{\partial\Theta^d}\,h_\rho(\Upsilon_\Theta x)\right|_{\Theta=0},
$$
$x\in\mathcal S$. Then
$$
G_\rho^{[\mathfrak d]}=\rho^\mathfrak d\sum_{d=0}^\mathfrak d\gamma_d^\mathfrak d\,g_\rho^{(d)}
$$
is called the modified directional Poisson wavelet of order~$\mathfrak d\in\mathbb N_0$ and
$$
H_\rho^{[\mathfrak d]}=\sum_{d=0}^\mathfrak d\gamma_d^\mathfrak d\,h_\rho^{(d)}
$$
is called inversion (or reconstruction) wavelet to the modified directional Poisson wavelet of order~$\mathfrak d\in\mathbb N_0$ if $\sum_{d=0}^\mathfrak d\gamma_d^\mathfrak d\,g_\rho^{(d)}$ and  $H_\rho^{[\mathfrak d]}$ satisfy~\eqref{eq:thesis_conjecture}.
\end{df}

\begin{lem}For each $d\in\mathbb N_0$ and each $\rho\in\mathbb R_+$ $h_\rho^{(d)}$ is an $\mathcal L^2(\mathcal S^n)$--function.\end{lem}
\begin{bfseries}Proof.\end{bfseries} Similar to the proof of Proposition~\ref{prop:well_definiteness}.\hfill$\Box$

\begin{thm} For a fixed~$\mathfrak d$, the families $\{G_\rho^{[\mathfrak d]}\}$ and $\{C\cdot H_\rho^{[\mathfrak d]}\}$ with the constant $C=\frac{\Sigma_n^2}{(n-1)^\mathfrak d\Gamma(\mathfrak d)}$, where~$G_\rho^{[\mathfrak d]}$ and $H_\rho^{[\mathfrak d]}$ are the modified directional Poisson wavelet of order~$\mathfrak d$ and its inversion wavelet, are an admissible wavelet pair with respect to $\alpha(\rho)=\frac{1}{\rho}$.
\end{thm}

\begin{bfseries}Proof. \end{bfseries}
According to~\eqref{eq:thesis_conjecture},
\begin{align}
\sum_{\kappa=1}^{N(n,l)}&\int_0^\infty\overline{a_l^\kappa(G_\rho^{[\mathfrak d]})}\,a_l^\kappa(H_\rho^{[\mathfrak d]})\,\frac{d\rho}{\rho}
   =\int_0^\infty\sum_{k_1=1}^{l}\overline{a_l^{k_1}(G_\rho^{[\mathfrak d]})}\,a_l^{k_1}(H_\rho^{[\mathfrak d]})\,\frac{d\rho}{\rho}\notag\\
&=\int_0^\infty\rho^{\mathfrak d}\,\overline{a_l^0(g_\rho^{(0)})}\,a_l^0(h_\rho^{(0)})\,[l(2\lambda+l)]^\mathfrak d\,\frac{d\rho}{\rho}.\label{eq:int_of_coeffs}
\end{align}
Since  $g_\rho^{(0)}$ and $h_\rho^{(0)}$ are given by the series
$$
g_\rho^{(0)}=\frac{1}{\Sigma_n}\sum_{l=0}^\infty\frac{\lambda+l}{\lambda}\,\frac{e^{-\rho l}}{A_l^0}\,Y_l^0,\qquad
h_\rho^{(0)}=\frac{1}{\Sigma_n}\sum_{l=0}^\infty\frac{\lambda+l}{\lambda}\,\frac{e^{-\frac{\rho l^2}{2\lambda}}}{A_l^0}\,Y_l^0,
$$
expression~\eqref{eq:int_of_coeffs} is equal to
\begin{align*}
&\frac{1}{\Sigma_n^2}\overbrace{\left(\frac{\lambda+l}{\lambda A_l^0}\right)^2}^{=N(n,l)}[l(2\lambda+l)]^\mathfrak d
   \int_0^\infty\rho^{\mathfrak d}\,e^{-\frac{\rho l(2\lambda+l)}{2\lambda}}\,\frac{d\rho}{\rho}\\
&\qquad=\frac{N(n,l)}{\Sigma_n^2}\,[l(2\lambda+l)]^\mathfrak d\left[\frac{2\lambda}{l(2\lambda+l)}\right]^\mathfrak d\Gamma(\mathfrak d)=\frac{N(n,l)}{C},
\end{align*}
i.e., $\{G_\rho^{[\mathfrak d]}\}$ and $\{C\cdot H_\rho^{[\mathfrak d]}\}$ satisfy condition~\eqref{eq:admbwv1}. In  order to verify condition~\eqref{eq:admbwv2}, we compute $\overline{G_\rho^{[\mathfrak d]}}\hat\ast H_\rho^{[\mathfrak d]}$ according to~\eqref{eq:zonal_product_Fc},
$$
\overline{G_\rho^{[\mathfrak d]}}\hat\ast H_\rho^{[\mathfrak d]}
   =\sum_{l=0}^\infty\sum_{k_1=1}^l\frac{\overline{a_l^{k_1}(G_\rho^{[\mathfrak d]})}\,a_l^{k_1}(H_\rho^{[\mathfrak d]})}{N(n,l)}\,\mathcal K_l,
$$
and further, by~\eqref{eq:thesis_conjecture},
\begin{align*}
\overline{G_\rho^{[\mathfrak d]}}\hat\ast H_\rho^{[\mathfrak d]}
   &=\rho^\mathfrak d\sum_{l=0}^\infty\frac{1}{N(n,l)}\,\overline{a_l^0(g_\rho^{(0)})}\,a_l^0(h_\rho^{(0)})\,[l(2\lambda+l)]^\mathfrak d\,\mathcal K_l\\
&=\frac{\rho^\mathfrak d}{\Sigma_n^2}\sum_{l=0}^\infty\underbrace{\frac{1}{N(n,l)}\left(\frac{\lambda+l}{\lambda\,A_l^0}\right)^2}_{=1}
   e^{-\frac{\rho l(2\lambda+l)}{2\lambda}}\,[l(2\lambda+l)]^\mathfrak d\,\mathcal K_l
\end{align*}
Since $|\mathcal K_l(t)|\leq\mathcal K_l(1)=\mathcal O(l^{2\lambda-1})$ for $l\to\infty$, cf. \cite[formula~8.937.4]{GR}, the series converges absolutely, and the order of summation and integration can be changed, i.e.
\begin{align*}
\int_R^\infty\overline{G_\rho^{[\mathfrak d]}}\hat\ast H_\rho^{[\mathfrak d]}\,\frac{d\rho}{\rho}&=\frac{1}{\Sigma_n^2}\sum_{l=0}^\infty
   \int_R^\infty[\rho l(2\lambda+l)]^\mathfrak d\,e^{-\frac{\rho l(2\lambda+l)}{2\lambda}}\,\frac{d\rho}{\rho}\cdot\mathcal K_l\\
&=\frac{1}{\Sigma_n^2}\sum_{l=0}^\infty W_{\mathfrak d-1}(Rl(2\lambda+l))\,e^{-\frac{R l(2\lambda+l)}{2\lambda}}\,\mathcal K_l,
\end{align*}
where~$W_{\mathfrak d-1}$ is a polynomial of degree~$\mathfrak d-1$. Consequently, it is a linear combination of Mexican needlets $K_{\sqrt R}^r$, $r=1,2,\dots,\mathfrak d-1$, cf.~\cite{aM10}, and the kernel~$K_{\sqrt R}^0$ of the operator $f(\rho\Delta^{\!\ast})$ for $f(s)=e^{-s}$. According to~\cite[Lemma~4.1 and Remark~4.2]{GM09a}, for every integer $N\geq0$ there exists $C_{r,N}$ such that
$$
R^{n/2}\left|\left(\frac{\theta}{\sqrt R}\right)^NK_{\sqrt R}^r(\cos\theta)\right|\leq C_{r,N}
$$
for all $\theta\in[0,\pi]$ and $R>0$ (in the case of a positive~$r$), respectively $0<R\leq1$ (for $r=0$). Consequently,
\begin{equation}\label{eq:estimation_zonal_product}
\left|\int_R^\infty\overline{G_\rho^{[\mathfrak d]}}\hat\ast H_\rho^{[\mathfrak d]}(\cos\theta)\,\frac{d\rho}{\rho}\right|\leq\frac{\mathfrak c R^{\frac{N-n}{2}}}{\theta^N}
\end{equation}
for some positive constant~$\mathfrak c$, every $N\in\mathbb N$, and $R\in(0,1]$. Now, in order to integrate the left--hand--side of~\eqref{eq:estimation_zonal_product} over~$\mathcal S^n$, divide the integration region into three parts: $0\leq\theta<R^{1/2}$ (integral~$I_1$), $R^{1/2}\leq\theta<R^{-1}$ (integral~$I_2$), and $R^{-1}\leq\theta\leq\pi$ (integral~$I_3$). For $r\leq1$ we obtain the following estimations:
\begin{align*}
I_1&\leq c\,R^{-n/2}\int_0^{R^{1/2}}\sin^{n-1}\theta\,d\theta\leq c\,R^{-n/2}\int_0^{R^{1/2}}\theta^{n-1}\,d\theta\leq c&(N=0),\\
I_2&\leq c\,R^{n/2}\int_{R^{1/2}}^{R^{-1}}\frac{\sin^{n-1}\theta}{\theta^{2n}}\,d\theta
   \leq c\,R^{n/2}\int_{R^{1/2}}^{R^{-1}}\theta^{-n-1}\,d\theta\leq c R^{3n/2}+c&(N=2n),\\
I_3&\leq c\int_{R^{-1}}^\pi\frac{\sin^{n-1}\theta}{\theta^n}\,d\theta\leq c\int_{R^{-1}}^\pi\frac{d\theta}{\theta^n}
   \leq c+c\,R^{n-1}&(N=n),
\end{align*}
where $c=$const. Therefore, \eqref{eq:admbwv2} holds for $r\in(0,1]$.

For $R>1$ each Mexican needlet is bounded by a constant, in order to see that, substitute $N=0$ to~\eqref{eq:estimation_zonal_product}. Further, for $r=0$ we have
\begin{align*}
\left|\sum_{l=0}^\infty e^{-\frac{Rl(2\lambda+l)}{2\lambda}}\mathcal K_l(t)\right|&\leq\left|\sum_{l=0}^\infty e^{-\frac{Rl(2\lambda+l)}{2\lambda}}\mathcal K_l(1)\right|\\
&\leq\left|\sum_{l=0}^\infty Rl(l(2\lambda+l))\,e^{-\frac{Rl(2\lambda+l)}{2\lambda}}\mathcal K_l(1)\right|\leq c.
\end{align*}
Consequently, the left--hand--side of~\eqref{eq:estimation_zonal_product} is bounded by a constant, hence, it is integrable over~$\mathcal S^n$ with integral value independent of~$R$.\hfill$\Box$

\begin{bfseries}Remark.\end{bfseries} In~\cite{HH09} the wavelets arise more naturally as directional derivatives of Poisson kernel (multiplied by a power of~$\rho$). The reason why we introduce linear combinations of those derivatives is the strong condition~\eqref{eq:admbwv1} which we want to be satisfied. Not only the simplicity gets lost in this case. Hayn and Holschneider show that directional derivatives of Poisson kernel have a good frequency localization, cf. discussion in \cite[Section~IV.B]{HH09}. We abandon these properties in order to obtain the exact reconstruction, contrary to the reconstruction with Fourier muliplier described in \cite[Section~I]{HH09}. For the same reason a wavelet family different from the Poisson wavelet family is needed for the inversion.

\section{Euclidean limit}\label{sec:Euclidean_limit}

A very important feature of directional Poisson wavelets is the existence of their Euclidean limit, i.e., the fact that for small scales the wavelets behave like wavelets over the Euclidean space. More exactly, for each~$d\in\mathbb N$ a square integrable function $G^{[d]}:\,\mathbb R^n\to\mathbb C$ exists such that
\begin{equation}\label{eq:Euclidean_limit}
\lim_{\rho\to0}\rho^ng_\rho^{[d]}\left(S^{-1}(\rho\xi)\right)=G^{[d]}(\xi)
\end{equation}
holds point-wise for every $\xi\in\mathbb R^n$, where~$S^{-1}$ is the inverse stereographic projection. Sufficient conditions for a wavelet family to have this property are stated in \cite[Theorem~3.4]{IIN14CWT} and verified for two--dimensional directional Poisson wavelets in Example following the proof of that theorem. In this section we want to compute Euclidean limits of directional Poisson wavelets.

\begin{thm}\label{thm:Euclidean_limit}
The Euclidean limits of directional Poisson wavelets are given by
\begin{equation}\label{eq:Euclidean_limit_PW}
G^{[d]}(\xi)=\frac{\partial^d}{\partial\xi_2^d}\,\frac{2}{\Sigma_n(1+|\xi|^2)^{\lambda+1}},
\end{equation}
$\xi=(\xi_2,\xi_3,\dots,\xi_{n+1})\in\mathbb R^n$.
\end{thm}

\begin{bfseries}Proof.\end{bfseries} According to~\cite[Theorem~7.1]{IIN14PW},
$$
G^{[0]}(\xi)=\frac{2}{(1+|\xi|^2)^{\lambda+1}}.
$$
Similarly as in the proof of \cite[Theorem~3.4]{IIN14CWT}, write $\xi\in\mathbb R^n$ in spherical coordinates,
$$
\xi=(R,\theta_2,\theta_3,\dots,\theta_{n-1},\varphi)
$$
with $R=|\xi|$ and $(\theta_2,\theta_3,\dots,\theta_{n-1},\varphi)$ -- spherical coordinates of an $n-1$-- dimensional sphere of radius~$r$. Denote by~$\theta$ the $\theta_1$--coordinate of $x=S^{-1}(\rho\xi)$, then
$$
\theta=2\arctan\frac{\rho R}{2}=\rho R+\mathcal O(\rho^3),\qquad\rho\to0,
$$
and further,
\begin{align*}
\sin\theta&=\sin(\rho R)+\mathcal O(\rho^3),\\
\cos\theta&=\cos(\rho R)+\mathcal O(\rho^4).
\end{align*}
for $\rho\to0$. Consequently,
\begin{align*}
g_\rho^{[d+1]}&\left(\Upsilon_\Theta S^{-1}(\rho\xi)\right)=\left.\rho\,\frac{\partial}{\partial\Theta}\,g_\rho^{[d]}\left(\Upsilon_\Theta S^{-1}(\rho\xi)\right)\right|_{\Theta=0}\\
&=\left.\rho\,\frac{\partial}{\partial\Theta}\,g_\rho^{[d]}\,\left(
\left(\begin{array}{ccccc}\cos\Theta&-\sin\Theta&0&\cdots&0\\
   \sin\Theta&\cos\Theta&0&\cdots&0\\
   0&0&1&\cdots&0\\
   \vdots&\vdots&\vdots&\ddots&\vdots\\
   0&0&0&\cdots&1\end{array}\right)
\left(\begin{array}{c}\cos(\rho R)+\mathcal O(\rho^4)\\\left(\sin(\rho R)+\mathcal O(\rho^3)\right)\cos\theta_2\\\sin(\rho R)\sin\theta_2\cos\theta_3\\\vdots\\
   \sin(\rho R)\sin\theta_2\dots\sin\theta_{n-1}\sin\varphi\end{array}\right)
   \right)\right|_{\Theta=0}\\
&=\left.\rho\,\frac{\partial}{\partial\Theta}\,g_\rho^{[d]}\,\left(
   \begin{array}{c}\cos\Theta\cos(\rho R)-\sin\Theta\sin(\rho R)\cos\theta_2+\mathcal O(\rho^4)\\
       \sin\Theta\cos(\rho R)+\cos\Theta\sin(\rho R)\cos\theta_2+\mathcal O(\rho^3)\\
       \sin(\rho R)\sin\theta_2\cos\theta_3+\mathcal O(\rho^3)\\\vdots\\\sin(\rho R)\sin\theta_2\dots\sin\theta_{n-1}\sin\varphi+\mathcal O(\rho^3)
   \end{array}\right)\right|_{\Theta=0}\\
&=\rho\left[\frac{\partial}{\partial x_1}\,g_\rho^{[d]}(\dots)[-\sin\Theta\cos(\rho R)-\cos\Theta\sin(\rho R)\cos\theta_2+\mathcal O(\rho^4)]\right.\\
&\qquad+\left.\frac{\partial}{\partial x_2}\,g_\rho^{[d]}(\dots)[\cos\Theta\cos(\rho R)-\sin\Theta\sin(\rho R)\cos\theta_2+\mathcal O(\rho^3)]\right]_{\Theta=0}\\
&=\rho\left[\frac{\partial}{\partial x_2}\,g_\rho^{[d]}(x)\cos(\rho R)+\mathcal O(\rho)\right]
\end{align*}
with
$$
x=\left(\begin{array}{c}\cos(\rho R)+\mathcal O(\rho^4)\\\sin(\rho R)\cos\theta_2+\mathcal O(\rho^3)\\\sin(\rho R)\sin\theta_2\cos\theta_3+\mathcal O(\rho^3)\\
   \vdots\\\sin(\rho R)\sin\theta_2\dots\sin\theta_{n-1}\sin\varphi+\mathcal O(\rho^3)\end{array}\right).
$$
Now, since $x_2=\sin(\rho R)\cos\theta_2=\rho\xi_2+\mathcal O(\rho^3)$, for $g_\rho^{[d+1]}$ we can write
\begin{align*}
g_\rho^{[d+1]}\left(\Upsilon_\Theta S^{-1}(\rho\xi)\right)&=\frac{\partial}{\partial\xi_2}\,g_\rho^{[d]}\left(\Upsilon_\Theta S^{-1}(\rho\xi)\right)\cos(\rho R)+\mathcal O(\rho^2)\\
&=\frac{\partial^d}{\partial\xi_2^d}\,p_\zeta\left(\Upsilon_\Theta S^{-1}(\rho\xi)\right)\cos(\rho R)+\mathcal O(\rho^2)
\end{align*}
for $\rho\to0$. Poisson kernels~$p_\zeta$, $\zeta=e^{-\rho}\hat e$, $\rho\in\mathbb R_+$, and their limit~$G^{[0]}$ are $\mathcal C^\infty$ functions, consequently, we can change the order of differentiation and limit calculation in~\eqref{eq:Euclidean_limit}. Thus,
\begin{align*}
G^{[d]}(\xi)&=\lim_{\rho\to0}\left.\rho^d\,\frac{\partial^d}{\partial\Theta^d}p_\zeta(\Upsilon_\Theta S^{-1}(\rho\xi))\right|_{\Theta=0}\\
&=\frac{\partial^d}{\partial\xi_2^d}\lim_{\rho\to0}p_\zeta(S^{-1}(\rho\xi))=\frac{\partial^d}{\partial\xi_2^d}G^{[0]}(\xi).
\end{align*}
\hfill$\Box$\\
\begin{bfseries}Remark.\end{bfseries} The mistaken factor~$n!$ in the expression representing the Euclidean limits of the directional Poisson wavelets in \cite[Theorem~IV.2]{HH09} occurs in the last equation of the proof on page~073512 in~\cite{HH09}: the equality is valid without~$n!$ on its right-hand-side.

\begin{exa} It is straightforward to compute the Euclidean limit of~$g_\rho^{[1]}$ for integer~$\lambda$ using~\eqref{eq:g_rho_1_explicit} and the l'Hospital rule. We have checked that the results coincide with~\eqref{eq:Euclidean_limit_PW} for $\lambda=\frac{1}{2}$, $\lambda=1$, $\lambda=\frac{3}{2}$, $\lambda=2$, $\lambda=3$, $\lambda=\frac{5}{2}$. In the case of non-integer~$\lambda$ we have computed the limit of~$(g_\rho^{[1]})^2$ and taken its negative square root as a result.
\end{exa}

\section*{Appendix}

Using recursive formulae for Fourier coefficients of directional derivatives of an $\mathcal L^2$--function given in Theorem~\ref{thm:series_representation} we are able to find an explicit expression for~$g_\rho^{[2]}$. Similarly as in Corollary~\ref{cor:f(1)explicit} let us compute the second directional derivative of a zonal function over~$\mathcal S^n$, $n\geq3$,
$$
f(x)=\sum_{l=0}^\infty\widehat f(l)\,C_l^\lambda(\cos\theta_1)=\sum_{l=0}^\infty\frac{\widehat f(l)}{A_l^0}\,Y_l^0(x),
$$
where $x=(\theta_1,\theta_2,\dots,\theta_{n-1},\phi)$ in spherical coordinates,
\begin{align}
f^{(2)}(x)&=-\sum_{l=0}^\infty\beta_{l,0}^2\,\frac{\widehat f(l)}{A_l^0}\,Y_l^0(x)+\sum_{l=2}^\infty\beta_{l,0}\beta_{l,1}\,\frac{\widehat f(l)}{A_l^0}\,Y_l^2(x)\notag\\
&=\frac{1}{2\lambda+1}\left[-\sum_{l=1}^\infty l(2\lambda+l)\,\widehat f(l)\,C_l^\lambda(\cos\theta_1)\right.\label{eq:2nd_der}\\
&\qquad\left.+\frac{8\lambda(\lambda+1)}{2\lambda-1}\sum_{l=1}^\infty\widehat f(l)\,C_{l-2}^{\lambda+2}(\cos\theta_1)\sin^2\theta_1
   \cdot C_2^{\lambda-\frac{1}{2}}(\cos\theta_2)\right],\notag
\end{align}
compare~\eqref{eq:beta_l_k1}, \eqref{eq:Alk1}, and~\eqref{eq:spherical_harmonics}. In order to simplify the calculations, we adapt the convention $C_l(t)=0$ for negative~$l$. Recursive formulae for Gegenbauer polynomials used in the sequel are valid also if negative indices occur. According to \cite[formula~8.933.3]{GR} we have
\begin{align*}
l(2\lambda+l)\,C_l^\lambda(t)&=2\lambda l\left[C_l^{\lambda+1}(t)-t\,C_{l-1}^{\lambda+1}(t)\right]\\
&=2\lambda\left[l\,C_l^{\lambda+1}(t)-(l-1)t\,C_{l-1}^{\lambda+1}(t)-t\,C_{l-1}^{\lambda+1}(t)\right],
\end{align*}
and, thus, by \cite[formula~8.933.2]{GR},
\begin{align*}
&l(2\lambda+l)\,C_l^\lambda(t)\\
&=4\lambda(\lambda+1)\left[t\,C_{l-1}^{\lambda+2}(t)-C_{l-2}^{\lambda+2}(t)
   -t^2\,C_{l-2}^{\lambda+2}(t)+t\,C_{l-3}^{\lambda+2}(t)\right]-2\lambda t\,C_{l-1}^{\lambda+1}(t)\\
&=4\lambda(\lambda+1)\left[t\left(C_{l-1}^{\lambda+2}(t)+C_{l-3}^{\lambda+2}(t)\right)-(1+t^2)\,C_{l-2}^{\lambda+2}(t)\right]-2\lambda t\,C_{l-1}^{\lambda+1}(t).
\end{align*}
Now, for~$C_2^{\lambda-\frac{1}{2}}$ we write
$$
C_2^{\lambda-\frac{1}{2}}(t)=\frac{1}{2}(2\lambda-1)\left[(2\lambda+1)\,t^2-1\right],
$$
cf.~\eqref{eq:C_jako_szereg} or~\cite[formula~8.930.3]{GR}, hence, the second summand in brackets on the right-hand-side of~\eqref{eq:2nd_der} is equal to
$$
4\lambda(\lambda+1)\sum_{l=1}^\infty\widehat f(l)\,C_{l-2}^{\lambda+2}(\cos\theta_1)\sin^2\theta_1\left[(2\lambda+1)\cos^2\theta_2-1\right].
$$
Consequently, $f^{(2)}(x)$ can be written as
\begin{align*}
f^{(2)}&(x)=\frac{1}{2\lambda+1}\biggl[-4\lambda(\lambda+1)\,\cos\theta_1
   \sum_{l=1}^\infty\widehat f(l)\left[C_{l-1}^{\lambda+2}(\cos\theta_1)+C_{l-3}^{\lambda+2}(\cos\theta_1)\right]\\
&+4\lambda(\lambda+1)(1+\cos^2\theta_1-\sin^2\theta_1)\sum_{l=1}^\infty\widehat f(l)\,C_{l-2}^{\lambda+2}(\cos\theta_1)\\
&+2\lambda\cos\theta_1\sum_{l=1}^\infty\widehat f(l)\,C_{l-1}^{\lambda+1}(\cos\theta_1)\biggr]\\
&+4\lambda(\lambda+1)\sum_{l=1}^\infty\widehat f(l)\,C_{l-2}^{\lambda+2}(\cos\theta_1)\sin^2\theta_1\cos^2\theta_2
\end{align*}
This expression can be simplified if we substitute $\frac{l+\lambda}{\Sigma_n\lambda}\,e^{-\rho l}$ for $\widehat f(l)$, i.e., compute the second directional derivative of Poisson kernel~$p_\zeta^\lambda$ located in $\zeta=e^{-\rho}\widehat e$ inside the unit sphere $\mathcal S^{2\lambda+1}$,
\begin{align*}
\Sigma_n&\left(p_\zeta^\lambda\right)^{(2)}(x)=\frac{1}{2\lambda+1}\\
&\cdot\left[-4\lambda(\lambda+1)\,\cos\theta_1
   \sum_{l=1}^\infty\frac{l+\lambda}{\lambda}\,e^{-\rho l}\left[C_{l-1}^{\lambda+2}(\cos\theta_1)+C_{l-3}^{\lambda+2}(\cos\theta_1)\right]\right.\\
&+8\lambda(\lambda+1)\cos^2\theta_1\cdot\frac{\lambda+2}{\lambda}\,e^{-2\rho}
   \sum_{l=1}^\infty\frac{(l-2)+(\lambda+2)}{\lambda+2}\,e^{-\rho(l-2)}C_{l-2}^{\lambda+2}(\cos\theta_1)\\
&\left.+2\lambda\cos\theta_1\cdot\frac{\lambda+1}{\lambda}\,e^{-\rho}
   \sum_{l=1}^\infty\frac{(l-1)+(\lambda+1)}{\lambda+1}\,e^{-\rho(l-1)}C_{l-1}^{\lambda+1}(\cos\theta_1)\right]\\
&+4\lambda(\lambda+1)\sin^2\theta_1\cos^2\theta_2\cdot\frac{\lambda+2}{\lambda}\,e^{-2\rho}\\
&\qquad\cdot\sum_{l=1}^\infty\frac{(l-2)+(\lambda+2)}{\lambda+2}\,e^{-\rho(l-2)}C_{l-2}^{\lambda+2}(\cos\theta_1)
\end{align*}
Now, consider the first series in this representation,
\begin{align*}
\sum_{l=1}^\infty&e^{-\rho l}\cdot\left[\frac{(l-1)+(\lambda+2)-1}{\lambda}\,C_{l-1}^{\lambda+2}(\cos\theta_1)
   +\frac{(l-3)+(\lambda+2)+1}{\lambda}\,C_{l-3}^{\lambda+2}(\cos\theta_1)\right]\\
&=\frac{\lambda+2}{\lambda}\,e^{-\rho}\sum_{l=1}^\infty \frac{(l-1)+(\lambda+2)}{\lambda+2}\,e^{-\rho(l-1)}\,C_{l-1}^{\lambda+2}(\cos\theta_1)\\
&+\frac{\lambda+2}{\lambda}\,e^{-3\rho}\sum_{l=1}^\infty\frac{(l-3)+(\lambda+2)}{\lambda+2}\,e^{-\rho(l-3)}\,C_{l-3}^{\lambda+2}(\cos\theta_1)\\
&-\sum_{l=1}^\infty\frac{e^{-\rho l}}{\lambda}\left[C_{l-1}^{\lambda+2}(\cos\theta_1)-C_{l-3}^{\lambda+2}(\cos\theta_1)\right],
\end{align*}
and use formula \cite[8.939.6]{GR} to compute $C_{l-1}^{\lambda+2}-C_{l-3}^{\lambda+2}$,
\begin{align*}
C_{l-1}^{\lambda+2}&(t)-C_{l-3}^{\lambda+2}(t)=C_{l-1}^{\lambda+2}(1)-C_{l-3}^{\lambda+2}(1)-2(l-2+\lambda+2)\,\int_t^1C_{l-2}^{\lambda+2}(\tilde t)\,d\tilde t.
\end{align*}
It yields
\begin{align*}
-\sum_{l=1}^\infty&\frac{e^{-\rho l}}{\lambda}\left[C_{l-1}^{\lambda+2}(\cos\theta_1)-C_{l-3}^{\lambda+2}(\cos\theta_1)\right]\\
&=\frac{e^{-3\rho}}{\lambda}\,\sum_{l=1}^\infty e^{-\rho(l-3)}\,C_{l-3}^{\lambda+2}(1)-\frac{e^{-\rho}}{\lambda}\,\sum_{l=1}^\infty e^{-\rho(l-1)}\,C_{l-1}^{\lambda+2}(1)\\
&+\frac{2(\lambda+2)\,e^{-2\rho}}{\lambda}\,\sum_{l=1}^\infty\frac{(l-2)+(\lambda+2)}{\lambda+2}\,e^{-\rho(l-2)}\int_{\cos\theta_1}^1C_{l-2}^{\lambda+2}(\tilde t)\,d\tilde t.
\end{align*}
In the first two summands of this expression the same series appears, and its value will be calculated with help of the generating function for Gegenbauer polynomials~\eqref{eq:generating_fct}. The last series converges absolutely, since the maximal value of~$C_{l-1}^{\lambda+2}$ increases polynomially in~$l$, and therefore the order of integration and summation can be changed. Consequently, for~$(p_\zeta^\lambda)^{(2)}$ we obtain the following expression
\begin{align*}
\left(p_\zeta^\lambda\right)^{(2)}&(x)=\frac{1}{2\lambda+1}\biggl[-4(\lambda+1)\,\cos\theta_1\biggl((\lambda+2)(e^{-\rho}+e^{-3\rho})\,p_\zeta^{\lambda+2}(x)\\
&+\frac{e^{-3\rho}-e^{-\rho}}{\Sigma_n}\cdot\frac{1}{(1-e^{-\rho})^{2(\lambda+2)}}
   +2(\lambda+2)\,e^{-2\rho}\int_{\cos\theta_1}^1p_\zeta^{\lambda+2}(\tilde t)\,d\tilde t\biggr)\\
&+8(\lambda+1)(\lambda+2)\cos^2\theta_1\,e^{-2\rho}\,p_\zeta^{\lambda+2}(x)+2(\lambda+1)\cos\theta_1\,e^{-\rho}\,p_\zeta^{\lambda+1}(x)\biggr]\\
&+4(\lambda+1)(\lambda+2)\sin^2\theta_1\cos^2\theta_2\,e^{-2\rho}\,p_\zeta^{\lambda+2}(x).
\end{align*}
Substituting~\eqref{eq:Poisson_kernel} yields
\begin{align*}
\left(p_\zeta^\lambda\right)^{(2)}(x)
   =&-\frac{2(\lambda+1)\,e^{-\rho}\,(1-e^{-2\rho})}{\Sigma_n\,(1-2e^{-\rho}\cos\theta_1+e^{-2\rho})^{\lambda+2}}\,\cos\theta_1\\
&+\frac{4(\lambda+1)(\lambda+2)\,e^{-2\rho}\,(1-e^{-2\rho})}{\Sigma_n\,(1-2e^{-\rho}\cos\theta_1+e^{-2\rho})^{\lambda+3}}\,\sin^2\theta_1\cos^2\theta_2
\end{align*}
Using this expression, one can compute the Euclidean limit of $g_\rho^{[2]}=\rho^2(p_\zeta^\lambda)^{(2)}$ for a fixed~$\lambda$, e.g.,
\begin{align*}
G^{[2]}(\xi)&=\frac{4\left(-1+2|\xi|^2+3|\xi|^2\cos(2\theta_2)\right)}{\pi^2\left(1+|\xi|^2\right)^4}&\text{for }\lambda=1,\\
G^{[2]}(\xi)&=\frac{12\left(-1+3|\xi|^2+4|\xi|^2\cos(2\theta_2)\right)}{\pi^3(1+|\xi|^2)^5}&\text{for }\lambda=2,\\
G^{[2]}(\xi)&=\frac{48\left(-1+4|\xi|^2+5|\xi|^2\cos(2\theta_2)\right)}{\pi^4(1+|\xi|^2)^6}&\text{for }\lambda=3.
\end{align*}
The results coincide with those obtained with formula~\eqref{eq:Euclidean_limit_PW}.

\end{document}